\documentclass[11pt]{amsart}
\usepackage{amssymb,mathtools,booktabs}
\usepackage[left=1.25in, right=1.25in, top=1.25in, bottom=1.25in]{geometry}
\usepackage[hidelinks]{hyperref}
\usepackage{aliascnt}
\usepackage[capitalise]{cleveref}
\crefname{appendix}{appendix}{appendices}
\Crefname{appendix}{Appendix}{Appendices}
\usepackage{todonotes}
\usepackage{tikz}
\newtheorem{theorem}{Theorem}[section]
\newaliascnt{proposition}{theorem}
\newtheorem{proposition}[proposition]{Proposition}
\aliascntresetthe{proposition}
\newaliascnt{lemma}{theorem}
\newtheorem{lemma}[lemma]{Lemma}
\aliascntresetthe{lemma}
\newaliascnt{corollary}{theorem}
\newtheorem{corollary}[corollary]{Corollary}
\aliascntresetthe{corollary}
\newaliascnt{conjecture}{theorem}
\newtheorem{conjecture}[conjecture]{Conjecture}
\aliascntresetthe{conjecture}
\theoremstyle{definition}
\newaliascnt{definition}{theorem}
\newtheorem{definition}[definition]{Definition}
\aliascntresetthe{definition}
\newaliascnt{example}{theorem}

\aliascntresetthe{example}
\theoremstyle{remark}
\newaliascnt{remark}{theorem}
\newtheorem{remark}[remark]{Remark}
\aliascntresetthe{remark}
\numberwithin{equation}{section}
\newcommand{\R}{\mathbb R}
\newcommand{\Z}{\mathbb Z}
\newcommand{\C}{\mathbb C}

\newcommand{\A}{\mathbb A}
\DeclareMathOperator{\Span}{span}
\DeclareMathOperator{\Sym}{Sym}
\DeclareMathOperator{\Frac}{Frac}

\DeclareMathOperator{\supp}{Supp}
\DeclareMathOperator{\Lie}{Lie}

\title{Reflection subgroups and the Lookup Conjecture}

\author{Christian Gaetz}
\address{Department of Mathematics, University of California, Berkeley, CA, USA.}
\email{\href{mailto:gaetz@berkeley.edu}{{\tt gaetz@berkeley.edu}}}
\author{Yibo Gao}
\address{Beijing International Center for Mathematical Research,\newline
Peking University, Beijing, China.}
\email{\href{mailto:gaoyibo@bicmr.pku.edu.cn}{{\tt gaoyibo@bicmr.pku.edu.cn}}}

\date{\today}

\begin{document}

\begin{abstract}
The \emph{Lookup Conjecture} of Boe--Graham (\emph{Amer. J. Math.}, 2003)
proposes a Bruhat graph criterion for rational smoothness of Schubert varieties in the flag variety $\mathcal{G}/\mathcal{B}$ of a complex reductive group. Its extension to Kac--Moody groups of affine type $\widetilde{A}_n$ was considered by Graham--Li. We prove the original conjecture as well as its extension to $\widetilde{A}_n$. Our method could be applied more generally to any Kac--Moody group after checking the rank-three reflection subgroups of its Weyl group.

We also prove a parabolic version of the Lookup Conjecture for the type $A$ affine Grassmannian, where a stronger Bruhat graph degree monotonicity property in fact holds, and discuss the implication via the geometric Satake correspondence for weight multiplicities.
\end{abstract}


\maketitle
\raggedbottom
\section{Introduction}

\subsection{The Lookup Conjecture} Let $W$ be the Weyl group of a Kac--Moody group $\mathcal{G}$, with length function $\ell$, Bruhat order $\le$, and reflection set $T$. For $x\le w$, write $I_w$ for the lower Bruhat interval $[e,w]$, $d_w(x)=\lvert\{t\in T \mid tx\le w\}\rvert$, and $q_w(x)=d_w(x)-\ell(w)$. Thus $d_w(x)$ equals the degree of $x$ in the Bruhat graph $\Gamma_w$ on $I_w$ and also the number of $\mathcal{T}$-invariant curves in the Schubert variety $X_w=\overline{\mathcal{B}w\mathcal{B}/\mathcal{B}}$ through the $\mathcal{T}$-fixed point $x\mathcal{B}$. Here $\mathcal{B}$ denotes a Borel subgroup of $\mathcal{G}$ containing the maximal torus $\mathcal{T}$. We call $q_w(x)$ the \emph{defect} of $x$ in $I_w$. Deodhar's conjecture \cite{Deodhar} that $q_w(x)\ge0$ was proved for arbitrary Coxeter groups by Dyer \cite{DyerDeodhar}.

The singularities of $X_w$ are of great interest because they are measured by the \emph{Kazhdan--Lusztig polynomials} $P_{x,w}(q)$ \cite{KL} and thereby are closely related to the Hecke algebra $\mathcal{H}_W$ and to the representation theory of $\mathfrak{g}=\Lie(\mathcal{G})$ \cite{Beilinson-Bernstein, Brylinski-Kashiwara}. In particular, $P_{x,w}(q)=1$ if and only if $X_w$ is rationally smooth at $x\mathcal{B}$. 

The Kazhdan--Lusztig polynomials are determined by a complicated recursion which obscures many of their properties (such as the non-negativity of their coefficients for general Coxeter groups \cite{Elias-Williamson}); the \emph{Combinatorial Invariance Conjecture}, due to Lusztig and to Dyer \cite{DyerThesis}, which asserts that $P_{x,w}$ depends only on the Bruhat graph on $[x,w]$, has received significant attention\footnote{A claimed proof of the conjecture was released by OpenAI in October 2026.}. One notable success in detecting the subtle geometry of $X_w$ directly in terms of $\Gamma_w$ is the celebrated \emph{Carrell--Peterson criterion} \cite{Carrell,Kumar}, which says that $X_w$ is rationally smooth at $x\mathcal{B}$ if and only if
\begin{equation}
\label{eq:CP-criterion}
q_w(y)=0 \qquad \text{for all $x \leq y \leq w$.}
\end{equation}

In type $A$, Deodhar proved the stronger result that $q_w(x)=0$ already implies rational smoothness at $x\mathcal{B}$ \cite{Deodhar}. This follows because in this case the tangent vectors to the $\mathcal{T}$-invariant curves counted by $d_w(x)$ form a basis for the tangent space to $X_w$ at $x\mathcal{B}$, and $\ell(w)=\dim(X_w)$. However, this property fails in other types, and the relationship between the invariant curves and full tangent space is still an active area of study \cite{Carrell-Kuttler-1, Carrell-Kuttler-2,Graham-Kreiman,Kisin-Pappas-Zhou}. The Lookup Conjecture proposes a much more parsimonious criterion than (\ref{eq:CP-criterion}) that it is hoped will characterize rational smoothness in general finite type.

\begin{conjecture}[Lookup Conjecture \cite{BG}]
\label{conj:lookup}
Let $\mathcal{G}$ be a complex reductive group with Weyl group $W$, and let $x \leq w$. If 
\begin{equation}
\label{eq:lookup-hypothesis}
q_w(y)=0 \qquad \text{for all $y \in \{x\} \cup \{tx \mid x<tx \leq w, t \in T\}$}
\end{equation}
then $X_w$ is rationally smooth at $x\mathcal{B}$.
\end{conjecture}

In short, rather than checking the defect of all elements of $[x,w]$, the conjecture says that it suffices to check at $x$ and at the upward Bruhat graph neighbors of $x$. In light of the Carrell--Peterson criterion, one could state the conjecture purely in terms of $\Gamma_w$, as the assertion that (\ref{eq:lookup-hypothesis}) implies (\ref{eq:CP-criterion}). So we can equally well say that the Lookup Conjecture holds for $\mathcal{G}$ or for its Weyl group $W$. Our first main result reduces the Lookup Conjecture for a Kac--Moody group $\mathcal{G}$ to a verification of this assertion for \emph{reflection subgroups} $H \leq W$ of rank three. 

\begin{theorem}\label{thm:main}
Let $W$ be a crystallographic Coxeter group. If the Lookup
Conjecture holds in every reflection subgroup $H \leq W$ of intrinsic
rank three, then it holds in $W$.
\end{theorem}

Boe and Graham verified the conjecture in finite Weyl groups of
small rank \cite[p.~319]{BG}, so we obtain a proof of \Cref{conj:lookup}. 

\begin{corollary}\label{cor:finite}
The Lookup Conjecture holds in every finite Weyl group; that is, \Cref{conj:lookup} is true.
\end{corollary}

Graham and Li extended the consideration of the Lookup Conjecture to infinite type, settling it for certain special Schubert varieties in affine type $\widetilde A_2$ \cite{GL}. The conjecture was later fully resolved for $\widetilde{A}_2$ by Boe and Graham \cite{BGAffine}. This result gives the remaining rank-three base case needed to apply \Cref{thm:main} in type $\widetilde{A}_n$, so we obtain:

\begin{corollary}\label{cor:affine}
The Lookup Conjecture holds for type $\widetilde{A}_n$.
\end{corollary}

\subsection{Degree monotonicity and the Parabolic Lookup Conjecture}

Our second main theorem concerns degree monotonicity in type $A$
affine Grassmannian quotients.
For $J\subseteq S$, let $W^J$ denote the set of minimal-length representatives
of right $W_J$-cosets.
For $w\in W^J$, the parabolic Bruhat graph $\Gamma_w^J$ has vertex set
$[e,w]\cap W^J$ and edges $x\to y$ if
$txW_J=yW_J$ for some $t\in T$ and $\ell(y)>\ell(x)$.
Write $d_w^J(x)$ for the degree of $x$ in $\Gamma_w^J$.

\begin{theorem}\label{thm:grassmannian}
Let $(W,S)$ have affine type $\widetilde A_n$, with $n\ge1$,
and let $J=S\setminus\{s\}$ for some $s\in S$.
For every $w\in W^J$,
\[
 d_w^J(x)\ge d_w^J(z)
 \qquad(x\le z\le w,\quad x,z\in W^J).
\]
\end{theorem}

\Cref{thm:grassmannian} implies the strengthened parabolic Lookup Conjecture of
Boe and Graham \cite[Conj.~1.2]{BG}.
When $W_J$ is finite, write $M_J$ for the set of maximal
representatives of right $W_J$-cosets. For $x\le w$ in $M_J$,
their assertion asks whether positive defect somewhere in $[x,w]$
forces positive defect at $x$ or at an upward reflection
neighbor that also belongs to $M_J$.
In the type $A$ affine Grassmannian, checking at $x$ suffices:

\begin{corollary}\label{cor:grassmannian-lookup}
In the setting of \Cref{thm:grassmannian}, for
$x\le w$ in $M_J$,
\[
 q_w(x)=0
 \quad\Longleftrightarrow\quad
 q_w(z)=0\text{ for every }x\le z\le w.
\]
In particular, the parabolic Lookup Conjecture holds for type $A$ affine Grassmannians.
\end{corollary}

The proofs of \Cref{thm:grassmannian} and \Cref{cor:grassmannian-lookup}
are independent of the reflection-subgroup reduction and
the rank-three computations and are special to type $\widetilde{A}_n$. The analog of \Cref{thm:grassmannian} does not hold in other affine types.

\subsubsection{Implications via Geometric Satake}
For \emph{spherical} Schubert varieties, the Geometric Satake Correspondence relates \Cref{thm:grassmannian} to weight multiplicities.
Take $\mathcal{G}=\operatorname{SL}_{n+1}(\C)$ and set
$\mathcal O=\C[[t]]$, $\mathcal K=\C((t))$, and let 
$\operatorname{Gr}_{\mathcal{G}}=\mathcal{G}(\mathcal K)/\mathcal{G}(\mathcal O)$ denote the affine Grassmannian.
For a dominant coweight $\lambda$ of $\mathcal{G}$, let
\[
 Y_\lambda
 =\overline{\mathcal{G}(\mathcal O)\tau^\lambda \mathcal{G}(\mathcal O)/\mathcal{G}(\mathcal O)}
 \subseteq\operatorname{Gr}_{\mathcal{G}},
\]
where $\tau^\lambda$ is the element of the affine Weyl group acting on $\mathfrak{h}_{\mathbb{R}}$ via translation by $\lambda$.
These spherical Schubert varieties are among the affine
Grassmannian Schubert varieties covered by
\Cref{thm:grassmannian}. Geometric Satake identifies the
category of $\mathcal{G}(\mathcal O)$-equivariant perverse sheaves on
$\operatorname{Gr}_{\mathcal{G}}$, with complex coefficients, with the
category of finite-dimensional representations of
$\mathcal{G}^\vee=\operatorname{PGL}_{n+1}(\C)$.
Under this equivalence, the intersection cohomology sheaf of
$Y_\lambda$ corresponds to the irreducible representation
$V_\lambda$ of highest weight $\lambda$
\cite[Thm.~7.3]{MV}. 

Write $W_J$ for the finite Weyl group of $\mathcal{G}$ and let $w_\nu$ denote the unique Bruhat-maximal element of $W_J\tau^\nu W_J$. In \Cref{cor:weight-count} we apply \Cref{cor:grassmannian-lookup} to show that
\[
 \dim V_\lambda[\mu]=1
 \quad\Longleftrightarrow\quad
 q_{w_\lambda}(w_\mu)=0.
\]
Thus a single Bruhat graph degree count can be used to detect one-dimensional weight spaces. \Cref{cor:weight-count} also contains more general implications for the weight
diagram of $V_\lambda$.

\subsubsection{Monotonicity in finite types}
In light of \Cref{thm:grassmannian}, we further classify all the parabolic quotients in finite Coxeter groups that are degree-monotone.
\begin{theorem}\label{thm:degree-monotone-finite-types}
Let $(W,S)$ be a finite irreducible Coxeter group. The following is the complete list of $W$ and $J \subsetneq S$ such that $d_w^J$ is monotone under Bruhat order for all $w \in W^J$:
\begin{enumerate}
\item $W$ is simply-laced; $S\setminus J$ contains only cominuscule nodes.
\item $W$ is of type $B_n$; $W_J$ is of type $B_{n-1}$ or $A_{n-1}$.
\item $(W,W_J)$ is $(H_3,I_2(5))$, or $W$ is of type $I_2(M)$ and $J$ is arbitrary. 
\end{enumerate}
\end{theorem}
Interestingly, the classification in (1) agrees, in the case of maximal parabolic quotients, with Stembridge's classification of \emph{fully commutative} quotients  \cite[Thm.~6.1]{stembridge-fullycommutative}. However \Cref{thm:grassmannian} implies that these notions diverge in affine type.

Our argument for the positive direction of \Cref{thm:degree-monotone-finite-types} is a stronger version of \cite[Thm.~1.4]{BG}, with a simpler and type-uniform proof. 

\subsection*{Outline} \Cref{sec:preliminaries} gives background and conventions on root systems and reflection subgroups and introduces \emph{slices}, the intersections of reflection subgroup cosets with lower Bruhat intervals. \Cref{sec:degree-comparison} applies the nil-Hecke ring of Kostant--Kumar to compare defects in $W$ with defects in reflection subgroups, and
\Cref{sec:geometry} uses Schubert geometry to control the maximal
components of a slice. \Cref{sec:reduction} combines these ideas to prove \Cref{thm:main} and  \Cref{cor:finite,cor:affine}. \Cref{sec:grassmannian} proves \Cref{thm:grassmannian} combinatorially,
then deduces \Cref{cor:grassmannian-lookup} and \Cref{cor:weight-count}. \Cref{thm:degree-monotone-finite-types} is proven in \Cref{sub:degree-monotone-finite-type}.

\section{Preliminaries}\label{sec:preliminaries}

We refer the reader to \cite{BB} for general background and notation for Coxeter groups. All Coxeter systems $(W,S)$ are assumed to be of finite rank, unless otherwise specified, with $S=\{s_1,\ldots,s_n\}$ denoting the set of simple generators and $W$ determined by relations $(s_is_j)^{m_{ij}}=e$ for a choice of integers $m_{ij} \in \{1,2,\ldots\} \cup \{\infty\}$ with $m_{ij}=m_{ji}$ for all $i,j$ and $m_{ij}=1$ if and only if $i=j$.

We write $x\lessdot y$ for a Bruhat cover and
$x\to y$ for an upward edge in the Bruhat graph $\Gamma$ on $W$: $y=tx>x$ for some $t\in T$. We use left reflection labels for edges unless otherwise specified. The \emph{support} $\supp(w)$ is the set of simple generators
in some (equivalently, every) reduced word for $w$. 

\subsection{Root systems}\label{sec:root-systems}
There are several distinct conventions for the root system of a Coxeter group, and we will make use of two.

\subsubsection{Geometric realizations}
A \emph{geometric realization} (in the sense of \cite[\S2]{DyerDeodhar}) of a Coxeter system $(W,S)$ is a realization of $W$ as a reflection group on a finite-dimensional real vector space $\tilde V$ with a nondegenerate $W$-invariant symmetric bilinear form $(,)$, together with a choice of linearly independent \emph{simple roots} $\alpha_1,\ldots,\alpha_n\in\tilde V$ satisfying $(\alpha_i,\alpha_i)=1$, $(\alpha_i,\alpha_j)=-\cos(\pi/m_{ij})$ if $m_{ij}<\infty$, and $(\alpha_i,\alpha_j)\le-1$ if $m_{ij}=\infty$. Each simple reflection $s_i$ acts by $v\mapsto v-2(v,\alpha_i)\alpha_i$. Every Coxeter system admits a geometric realization.

The \emph{roots} are the unit vectors $w\alpha_i$ for $w\in W$ and $1\le i\le n$. Each root lies in the nonnegative or the nonpositive span of the simple roots, and is called \emph{positive} or \emph{negative} accordingly. The reflection of a root $\alpha=w\alpha_i$ is $s_\alpha=ws_iw^{-1}$, which acts by $v\mapsto v-2(v,\alpha)\alpha$. 

\subsubsection{Kac--Moody realizations}
For an $n\times n$ integral (generalized) Cartan matrix $A=(a_{ij})$ of rank $r$, a \emph{Kac--Moody realization} is the real vector space $\mathfrak h_\R$ of dimension $2n-r$, together with linearly independent simple coroots $\alpha_i^\vee\in\mathfrak h_\R$ and simple roots $\alpha_i\in\mathfrak h_\R^*$ satisfying $\langle\alpha_i^\vee,\alpha_j\rangle=a_{ij}$ \cite[\S1.1]{Kumar}. The Weyl group acts on $\mathfrak h_\R^*$ by $s_i(\lambda)=\lambda-\langle\alpha_i^\vee,\lambda\rangle\alpha_i$, and $\mathfrak h_\R\otimes\C$ is the Lie algebra of a maximal torus of a corresponding Kac--Moody group $\mathcal{G}$. The Weyl group $W$ is a Coxeter group with $4\cos^2\left(\frac{\pi}{m_{ij}}\right)=a_{ij}a_{ji}$ for $m_{ij}$ finite and $a_{ij}a_{ji} \geq 4$ for $m_{ij}=\infty$. Thus $(W,S)$ admits a Kac--Moody realization if and only if all $m_{ij}$ lie in $\{1,2,3,4,6,\infty\}$; such $W$ are called \emph{crystallographic}. While many of our results work for arbitrary Coxeter groups, the algebrogeometric methods of \Cref{sec:geometry} require crystallography.

The (real) roots in this realization are $w\alpha_i\in\mathfrak h_\R^*$ for $w\in W$ and $1\le i\le n$; the imaginary roots of the Kac--Moody algebra play no role in this paper. Each root is an integral combination of simple roots whose coefficients are all nonnegative or all nonpositive, and is called \emph{positive} or \emph{negative} accordingly. A root $\alpha=w\alpha_i$ has corresponding coroot $\alpha^\vee=w\alpha_i^\vee\in\mathfrak h_\R$, which depends only on $\alpha$, and reflection $s_\alpha=ws_iw^{-1}$, which acts by $\lambda\mapsto\lambda-\langle\alpha^\vee,\lambda\rangle\alpha$.

\subsubsection{Notation} Having fixed a Coxeter system $(W,S)$ and either a geometric or Kac--Moody realization, we write $\Delta \subset \Phi^+ \subset \Phi$ for the simple roots, positive roots, and roots, respectively, and $\Phi^-=-\Phi^+$ for the negative roots. We write $V$ for the real span of $\Phi$ inside either $\tilde V$ or $\mathfrak{h}^*_\R$. Both geometric and Kac--Moody realizations of root systems are \emph{reduced}, meaning that $\R\alpha \cap \Phi = \{\pm \alpha\}$ for all roots $\alpha$.

\subsection{Reflection subgroups}\label{sec:intrinsic}

A subgroup $H$ of a Coxeter group $W$ is called a \emph{reflection subgroup} if it is generated by $H \cap T$. It was shown in \cite{DeodharSub} and \cite[Thm.~1.8]{DyerThesis} that $H$ is a Coxeter group with simple reflections
\begin{equation*}%\label{eq:intrinsic-generators}
 S_H \coloneqq \{t\in H\cap T \mid \ell(tt')>\ell(t) \text{ for all } t' \in (H\cap T) \setminus \{t\}\}
\end{equation*}
and reflections $T_H=H\cap T$. If $H \leq K \leq W$ are reflection subgroups, it follows by \emph{op. cit.} that $S_H$ does not depend on whether the ambient group is taken to be $K$ or $W$. We write $\ell_H,\le_H,$ and $\Gamma_H$ for the intrinsic length, Bruhat order, and Bruhat graph on $(H,S_H)$, respectively. We call $\lvert S_H\rvert$ the \emph{intrinsic rank} of $H$. Reflection subgroups of crystallographic Coxeter groups are crystallographic. 

For a subspace $U\subseteq V$ spanned by roots, write
\begin{equation}\label{eq:root-subgroup}
 H_U=\langle s_\alpha \mid \alpha\in\Phi\cap U\rangle
\end{equation}
for the corresponding reflection subgroup. Its reflection set is exactly
\begin{equation}\label{eq:root-reflections}
 T_{H_U}=\{s_\alpha \mid \alpha\in\Phi\cap U\}.
\end{equation}

\begin{remark}
    The reader is cautioned that the intrinsic rank of $H_U$ may in general differ from $\dim_\R(U)$. 
\end{remark}

\subsection{Slices of Bruhat intervals} The intersections of reflection group cosets with Bruhat order ideals will be of central importance in this paper. 

For a subset $X\subseteq W$ and $x \in X$, we write $N_X(x)$ for the set $\{tx\in X \mid t\in T\}$ of neighbors of $x$ in the Bruhat graph $\Gamma_X$ restricted to $X$, and $d_X(x)=\lvert N_X(x)\rvert$. 
We will most often be interested in the case $X=I$ is a Bruhat order ideal. We write $N_w$ and $\Gamma_w$ for $N_{I_w}$ and $\Gamma_{I_w}$, respectively, where $I_w=[e,w]$. We say that $I_w$ is \emph{regular above $x$} if $q_w(y)=0$ for
all $y\in[x,w]$; that is, if (\ref{eq:CP-criterion}) holds for $x$. We say that $x$ \emph{satisfies the Lookup
Hypotheses} (with respect to $w$) if (\ref{eq:lookup-hypothesis}) holds.

\begin{proposition}[\cite{DyerThesis,DyerGraph}]
\label{prop:Dyer-bg}
Let $H$ be any reflection subgroup of $W$ and set
$\Phi_H^+=\{\alpha\in\Phi^+ \mid s_\alpha\in H\}$. Any coset $Hx$ contains a unique minimal element $x_0$ under the Bruhat order on $W$ and this element is characterized by $x_0^{-1}(\Phi_H^+)\subseteq\Phi^+$. Furthermore, the map $h\mapsto hx_0$
identifies $\Gamma_H$ with $\Gamma_{Hx}$.
\end{proposition}

\begin{definition}\label{def:slice}
The \emph{slice} of $I_w$ on $Hx$ is
\[
I_H(w;x)\coloneqq \{h\in H \mid hx_0\le w\},
\]
where $x_0$ is as in \Cref{prop:Dyer-bg}. It is a finite order ideal $I$ for $\le_H$. For $u \in I$, we write
$d_I(u) \coloneqq \lvert\{t\in T_H \mid tu\in I\}\rvert$ and define
\[
 \mathcal M_I(u)=
 \{v\in\max_{\le_H}I \mid u\le_Hv\},
\]
the set of $\leq_H$-maximal elements of $I$ above $u$. For $v\in\mathcal M_I(u)$, the principal ideal $[e,v]_H$ is called
a \emph{maximal slice component through $u$}. We define a graph structure on $\mathcal M_I(u)$ with an edge between
distinct $v,v'$ when there is $h\in H$ such that
\[
 u\le_Hh,\qquad h\lessdot_Hv,\qquad h\lessdot_Hv'.
\]
\end{definition}

\section{Degree comparisons from nil-Hecke rings}\label{sec:degree-comparison}

In this section we apply the nil-Hecke ring of Kostant--Kumar
\cite{KostantKumar-86-nilhecke} to relate defects in $W$ to defects
in reflection subgroups. The application of this ring to the study of defects was pioneered in Dyer's proof of Deodhar's conjecture that defects are nonnegative \cite{DyerDeodhar}. Except where specified, $W$ is an arbitrary Coxeter group in a geometric realization; the conclusions of \Cref{prop:plane-comparison,prop:separation} involve only $(W,S)$, so they hold in every realization.

\subsection{Nil-Hecke coefficients}
We follow the conventions of \cite[\S2]{DyerDeodhar}. Let $W$ be a Coxeter group, acting in a geometric realization on $\tilde{V}$. Set $\mathcal{R}=\Sym(\tilde{V})$, graded with $\tilde{V}$ in degree one, and
$\mathcal{Q}=\Frac(\mathcal{R})$, with the induced left $W$-actions.

Let $\mathcal{Q}_W$ be the ring of \cite[(4.1)]{KostantKumar-86-nilhecke}: the free right $\mathcal{Q}$-module with
basis $\{\mathbf{b}_x \mid x\in W\}$ and multiplication
\[
 (\mathbf{b}_x F)(\mathbf{b}_y G)=\mathbf{b}_{xy}\,y^{-1}(F)G
 \qquad(F,G\in\mathcal{Q},\ x,y\in W).
\]
Identifying $F\in\mathcal{Q}$ with $\mathbf{b}_eF$ makes $\mathcal{Q}$ a subring of $\mathcal{Q}_W$, and
the multiplication rule gives
\begin{equation}
\label{eq:left-module}
 F\mathbf{b}_x=\mathbf{b}_x\,x^{-1}(F)\qquad(F\in\mathcal{Q},\ x\in W).
\end{equation}
Thus $\mathcal{Q}_W$ is also a free left $\mathcal{Q}$-module with basis
$\{\mathbf{b}_x \mid x\in W\}$. For $1\le i\le n$ define
\[
 \mathbf{x}_i=\alpha_i^{-1}(\mathbf{b}_{s_i}-\mathbf{b}_e).
\]
These elements satisfy the braid relations, so
$\mathbf{x}_w \coloneqq \mathbf{x}_{i_1}\cdots \mathbf{x}_{i_L}$ is independent of the
reduced word $w=s_{i_1}\cdots s_{i_L}$, with $\mathbf{x}_e=\mathbf{b}_e$.

Define $E_{x,w}\in\mathcal{Q}$ by
\[
 \mathbf{x}_w=\sum_{x \in W}E_{x,w}\mathbf{b}_x.
\]
By the subword property, $E_{x,w} =0$ if $x \not\leq w$, and, by
\cite[p.~573]{DyerDeodhar}, for $x\le w$ we have,
\begin{equation}\label{eq:nilhecke-numerator}
 E_{x,w}=\frac{C_{x,w}}{A_{x,w}},\qquad
 A_{x,w}=\prod_{\substack{\alpha\in\Phi^+\\s_\alpha x\le w}}\alpha,
 \qquad 0\ne C_{x,w}\in\mathcal{R},\quad \deg C_{x,w}=q_w(x),
\end{equation}
where $C_{x,w}$ is homogeneous.

For $\alpha \in \Phi$ and $F$ in the localization $\mathcal{R}_{(\alpha)}$ of $\mathcal{R}$ at the prime ideal $(\alpha)$, we write $\overline{F}$ for its image in $\mathcal{R}_{(\alpha)}/\alpha \mathcal{R}_{(\alpha)} \cong \Frac(\mathcal{R}/(\alpha))$. Since $A_{x,w}$ is a product of distinct
positive roots, $\alpha E_{x,w}$ lies in $\mathcal{R}_{(\alpha)} \subset \mathcal{Q}$.
We call $\overline{\alpha E_{x,w}}$ the \emph{residue} of $E_{x,w}$ along $\alpha$.

\begin{proposition}\label{prop:residue}
Let $x,w\in W$ and $\alpha\in\Phi$. Then
$E_{x,w}+E_{s_\alpha x,w}\in\mathcal{R}_{(\alpha)}$. Consequently,
\[
 \overline{\alpha E_{x,w}}+
 \overline{\alpha E_{s_\alpha x,w}}=0.
\]
If $x\le w$ and $q_w(x)=0$, then $C_{x,w}$ is a nonzero constant. If
moreover $\alpha\in\Phi^+$ and $s_\alpha x\le w$, then
$\overline{\alpha E_{x,w}}\ne0$.
\end{proposition}
\begin{proof}
The final two assertions are immediate from
\eqref{eq:nilhecke-numerator}: if $q_w(x)=0$, then $C_{x,w}$ is a
nonzero constant, and if moreover $\alpha\in\Phi^+$ and
$s_\alpha x\le w$, then $\overline{\alpha E_{x,w}}$ equals $C_{x,w}$
divided by the image of $A_{x,w}/\alpha$, a product of roots not
proportional to $\alpha$.

The second assertion follows from the first: both residues are
defined, and their sum is the image of
$\alpha(E_{x,w}+E_{s_\alpha x,w})\in\alpha\mathcal{R}_{(\alpha)}$.
We prove the first by induction on $\ell(w)$. If $w=e$, then
$\mathbf{x}_e=\mathbf{b}_e$, so each $E_{x,e}$ is $0$ or $1$.

Suppose $\ell(w)>0$, and let $s_i$ be a left descent of $w$. Then
$\mathbf{x}_w=\mathbf{x}_i\mathbf{x}_{w'}$ for $w'=s_iw$. We have $\mathbf{b}_{s_i}F=s_i(F)\mathbf{b}_{s_i}$ for
$F\in\mathcal{Q}$ by (\ref{eq:left-module}), so
\[
 \mathbf{x}_w
 =\alpha_i^{-1}(\mathbf{b}_{s_i}-\mathbf{b}_e)\sum_{y\in W}E_{y,w'}\mathbf{b}_y
 =\sum_{y\in W}\alpha_i^{-1}\bigl(s_i(E_{y,w'})\mathbf{b}_{s_iy}-E_{y,w'}\mathbf{b}_y\bigr).
\]
Comparing coefficients of $\mathbf{b}_x$ gives
\begin{equation}
    \label{eq:E-difference}
 E_{x,w}=\alpha_i^{-1}\bigl(s_i(E_{s_ix,w'})-E_{x,w'}\bigr)
 \qquad(x\in W).
\end{equation}
Let $\beta=s_i\alpha\in\Phi$. Then $s_\beta=s_is_\alpha s_i$, so
$s_is_\alpha x=s_\beta s_ix$. Adding (\ref{eq:E-difference}) at $x$ and
at $s_\alpha x$ gives
\[
 E_{x,w}+E_{s_\alpha x,w}=\alpha_i^{-1}\bigl(s_i(F')-F\bigr),
 \qquad
 F=E_{x,w'}+E_{s_\alpha x,w'},\quad
 F'=E_{s_ix,w'}+E_{s_\beta s_ix,w'}.
\]
By induction, $F\in\mathcal{R}_{(\alpha)}$ and
$F'\in\mathcal{R}_{(\beta)}$. The ring automorphism $s_i$ of
$\mathcal{R}$ carries the prime $(\beta)$ to $(\alpha)$, so
$s_i(F')\in\mathcal{R}_{(\alpha)}$.

If $\alpha\ne\pm\alpha_i$, then $\alpha_i$ is not proportional to
$\alpha$, so $\alpha_i^{-1}\in\mathcal{R}_{(\alpha)}$ and the claim
follows. If $\alpha=\pm\alpha_i$, then
$\mathcal{R}_{(\alpha)}=\mathcal{R}_{(\alpha_i)}$, $\beta=-\alpha$,
and $s_\beta=s_i$, so $F'=E_{s_ix,w'}+E_{x,w'}=F$. Since
$s_i(\alpha_i)=-\alpha_i$ and $s_i(v)-v=-2(v,\alpha_i)\alpha_i$ for
$v\in\tilde V$, the automorphism $s_i$ preserves the prime
$(\alpha_i)$ and induces the identity on $\mathcal{R}/(\alpha_i)$;
that is, $s_i(G)-G\in\alpha_i\mathcal{R}$ for $G\in\mathcal{R}$.
Writing $F=G/H$ with $G,H\in\mathcal{R}$ and $\alpha_i\nmid H$, we get
\[
 s_i(F)-F
 =\frac{\bigl(s_i(G)-G\bigr)H-G\bigl(s_i(H)-H\bigr)}{H\,s_i(H)},
\]
whose numerator lies in $\alpha_i\mathcal{R}$ and whose denominator
is not divisible by $\alpha_i$. Hence
$s_i(F)-F\in\alpha_i\mathcal{R}_{(\alpha_i)}$, and
$E_{x,w}+E_{s_\alpha x,w}=\alpha_i^{-1}\bigl(s_i(F)-F\bigr)$ lies in
$\mathcal{R}_{(\alpha_i)}=\mathcal{R}_{(\alpha)}$.
\end{proof}

\subsection{Contributions from dihedral reflection subgroups}\label{sec:dihedral-contributions}
We call a Coxeter group \emph{dihedral} if it has rank two. Given $y=s_{\alpha}x$, our approach to comparing $d_w(x)$ and $d_w(y)$ will be to compare the degrees of $x$ and $y$ within a coset of each \emph{maximal} dihedral reflection subgroup of $W$ containing $s_{\alpha}$. For a plane $P=\Span_\R\{\alpha,\beta\} \subset V$ with $\beta\in\Phi\setminus\R\alpha$, \cite[Lem.~3.17]{DyerThesis} implies that the reflection subgroup $H_P$ is dihedral. Moreover, \cite[Cor.~3.18]{DyerThesis} tells us that these are exactly the maximal dihedral reflection subgroups containing $s_{\alpha}$; this
description involves only $(W,S)$, so the collection $\mathcal{D}_\alpha$ of these subgroups $H_P$ does not depend on the root system realization. In either a geometric or a Kac--Moody realization, the roots of $M\in\mathcal{D}_\alpha$ lie in a single plane: $M$ is generated by two reflections $s_{\beta_1},s_{\beta_2}$, which preserve $\Span_\R\{\beta_1,\beta_2\}$, and every reflection of $M$ is conjugate in $M$ to one of these.

We first examine the case $W$ is dihedral; the proof of \Cref{prop:easy-dihedral} is easy and omitted.

\begin{proposition}
\label{prop:easy-dihedral}
Let $W$ be a dihedral Coxeter group. Then $x < y$ if and only if $\ell(x)<\ell(y)$, and $x \to y$ if and only if $\ell(y)-\ell(x)$ is positive and odd. A nonempty finite Bruhat order ideal $I$ is either principal with top length $k$, or is generated by both elements of length $k$. If $I=I_w$ is principal, we have $d_w(x)=k$ for all $x \in I$. Otherwise,
\[
d_I(x)=\begin{cases}
    k & \text{if $\ell(x)\equiv k \mod{2}$},\\
    k+1 & \text{otherwise.}
\end{cases}
\]
\end{proposition}

We now return to general Coxeter groups. For $M \in \mathcal{D}_\alpha$ and $z \in I_w$ define $d_{M,w}(z)\coloneqq \lvert\{t\in T_M \mid tz\le w\}\rvert$. For neighbors $x, y=s_\alpha x$ in $I_w$ let $\delta_M(x,y)=d_{M,w}(y)-d_{M,w}(x)$. The collection $\mathcal{D}_\alpha$ partitions the reflections other than $s_\alpha$, so
\begin{equation}\label{eq:plane-sum}
 d_w(y)-d_w(x)=\sum_{M\in \mathcal{D}_\alpha}\delta_M(x,y).
\end{equation}
Only finitely many terms in the sum are nonzero, since only finitely many reflections label edges of $\Gamma_w$.

\begin{proposition}\label{prop:plane-comparison}
Let $W$ be a Coxeter group in a geometric or Kac--Moody realization, and let
$x,w\in W$ and $\alpha\in\Phi^+$ satisfy $x,y=s_\alpha x\in I_w$. If $q_w(x)=0$, then $\delta_M(x,y)\in\{0,1\}$ for every
$M \in \mathcal{D}_\alpha$, and $q_w(y)=\sum_{M \in \mathcal{D}_\alpha}\delta_M(x,y)$.
In particular, if $q_w(y)=0$, then $\delta_M(x,y)=0$ for all $M$.
\end{proposition}
\begin{proof}
The statement depends only on $(W,S)$, so we may work
in a geometric realization. Each $M \in \mathcal{D}_\alpha$ equals $H_P$ for a root plane
$P\ni\alpha$, and distinct $M$ have distinct planes. Positive roots $\beta, \beta' \neq \alpha$ are proportional modulo $\alpha$ exactly when $\Span_\R\{\beta, \beta', \alpha\}$ is a plane, that is, when
$s_{\beta},s_{\beta'}$ lie in the same $M \in \mathcal{D}_\alpha$. For each $M$, fix $\beta_M\in\Phi^+\setminus\{\alpha\}$ with $s_{\beta_M}\in M$. By \Cref{prop:residue} and
\eqref{eq:nilhecke-numerator} we have
\[
 \overline{C_{y,w}}
 =-C_{x,w}\,\frac{\overline{(A_{y,w}/\alpha)}}{\overline{(A_{x,w}/\alpha)}}
 =c\prod_{M \in \mathcal{D}_\alpha}\overline{\beta_M}^{\,\delta_M(x,y)}
 \ne0,
\]
with $C_{x,w}$ and $c$ nonzero constants. Since $\overline{C_{y,w}} \in \mathcal{R}/(\alpha)$, we have $\delta_M(x,y)\ge0$ separately for each $M \in \mathcal{D}_\alpha$ by unique factorization.
The slice
$I_M(w;x)$ is a finite intrinsic order ideal whose Bruhat graph is isomorphic to $\Gamma_{Mx\cap I_w}$ by \Cref{prop:Dyer-bg}. In particular,
$d_{M,w}(x)$ and $d_{M,w}(y)$ are its degrees at the
endpoints of the edge corresponding to $x$--$y$, so
$\delta_M(x,y)$ is the degree difference along that edge. By \Cref{prop:easy-dihedral}, we therefore have $\lvert \delta_M(x,y) \rvert \leq 1$. Now apply
\eqref{eq:plane-sum}.
\end{proof}

\begin{corollary}\label{cor:internal}
Let $W$ be a Coxeter group in a geometric or Kac--Moody
realization. Let $w\in W$,
let $U\subseteq V$ be spanned by roots, set $H=H_U$, and
let $x_0$ be the minimum of a coset $Hx_0$ meeting $I_w$.
Let $I=I_H(w;x_0)$ and suppose $u\in I$ and $\alpha\in\Phi^+\cap U$ satisfy $u'=s_\alpha u\in I$.
Define $x=ux_0$ and $y=u'x_0$. If $q_w(x)=0$, then
\[
 d_I(u')-d_I(u)=\sum_{\substack{M \in \mathcal{D}_\alpha \\ M\subseteq H}}\delta_M(x,y)\ge0.
\]
Thus $d_I(u')-d_I(u)$ is $0$ if $q_w(y)=0$, and is positive
if some such $\delta_M(x,y)$ is nonzero.
\end{corollary}
\begin{proof}
By \Cref{def:slice} and \eqref{eq:root-reflections}, $d_I(u)$ (resp. $d_I(u')$) counts the
reflections $s_\gamma$ with $\gamma\in U$ and $s_\gamma x\le w$ (resp.
$s_\gamma y\le w$). If $M \in \mathcal{D}_\alpha$ contains some $s_\gamma\ne s_\alpha$
with $\gamma\in U$, then every root of $M$ lies in $\Span_\R\{\alpha,\gamma\}\subseteq U$,
so $M\subseteq H$. Hence the reflections of $H$ other than $s_\alpha$ are partitioned by
the subgroups $M\subseteq H$. The result follows by taking the corresponding sub-sum of \eqref{eq:plane-sum}
and using \Cref{prop:plane-comparison}.
\end{proof}

\subsection{Neighbor sets} We will also apply the nil-Hecke ring to prove \Cref{prop:separation} about sets of Bruhat graph neighbors. 

\begin{lemma}\label{lem:linear-section}
Suppose homogeneous linear polynomials $\{F_z\in\mathcal{R} \mid z\in[x,y]\}$ satisfy
$F_{z'}-F_z\in\R\alpha$ whenever $z'=s_\alpha z$ are adjacent in $\Gamma$. If $F_x=0$ and $F_{z'}=0$ at every cover
$x\lessdot z'\le y$, then $F_z=0$ for all $z$.
\end{lemma}
\begin{proof}
Induct on $\ell(x,z)$, the cases zero and one being by hypothesis.
If $\ell(x,z) \geq 2$, choose $z_0\in[x,z]$ with $\ell(z_0,z)=2$. 
Then there are two middle elements $z_1,z_2$ in the diamond $[z_0,z]$. Both $F_{z_1},F_{z_2}$ vanish
by induction. Thus $F_z$ lies in the two distinct root lines of
the edges $z_1 \to z$ and $z_2 \to z$, and so must be zero.
\end{proof}

\begin{lemma}\label{lem:coatom-ratio}
Suppose $x\le y\lessdot w$ and $I_w$ is regular above $x$.
For $z\in[x,y]$, the functions $F^w_z=E_{z,y}/E_{z,w}$
are nonzero homogeneous linear polynomials satisfying the
edge congruences of \Cref{lem:linear-section}. Moreover, $F^w_y$ is a nonzero multiple of the root of $y\to w$.
\end{lemma}
\begin{proof}
Since $I_y\subseteq I_w$, $A_{z,y}$ divides $A_{z,w}$, and $C_{z,w}$
is a nonzero constant by regularity and \Cref{prop:residue}. Hence
\[
 F^w_z=\frac{C_{z,y}}{C_{z,w}}\frac{A_{z,w}}{A_{z,y}}
\]
is polynomial. The functions $E_{z,y}$ and $E_{z,w}$ are
homogeneous of degrees $-\ell(y)$ and $-\ell(w)$, so $\deg F^w_z=-\ell(y)+\ell(w)=1$.

Consider an edge $z'=s_\alpha z$ in $[x,y]$. By \Cref{prop:residue}, the residues
$\overline{\alpha E_{z,w}},\overline{\alpha E_{z',w}}$ are nonzero and
\[
 \overline{F^w_z}=\overline{\alpha E_{z,y}}/\overline{\alpha E_{z,w}}=\overline{\alpha E_{z',y}}/\overline{\alpha E_{z',w}}=\overline{F^w_{z'}},
\]
which is the required congruence. Finally, $A_{y,y}$ and $A_{y,w}$ differ only in the factor coming from $y\to w$, and both $C_{y,y},C_{y,w}$ are constants by
\Cref{prop:residue}, proving the last assertion.
\end{proof}

\begin{proposition}\label{prop:separation}
Let $x\le y\lessdot w,w'$ with $w\ne w'$,
and suppose $I_w$ and $I_{w'}$ are regular above $x$. If
$N_w(x)=N_{w'}(x)$, then $N_w(z')\ne N_{w'}(z')$ for
some cover $x\lessdot z'\le y$.
\end{proposition}
\begin{proof}
Suppose instead that $N_w(x)=N_{w'}(x)$ and that
$N_w(z')=N_{w'}(z')$ for every cover $x\lessdot z'\le y$. Let $z$ be
$x$ or one of these covers. Then $A_{z,w}=A_{z,w'}$, since both
products in \eqref{eq:nilhecke-numerator} run over the
$\alpha\in\Phi^+$ with $s_\alpha z\in N_w(z)=N_{w'}(z)$, and
$C_{z,w}$ and $C_{z,w'}$ are nonzero constants by regularity and
\Cref{prop:residue}. Therefore
$E_{z,w'}=\lambda_zE_{z,w}$ for a nonzero scalar $\lambda_z$.
\Cref{prop:residue} applied to the edges $x \to z'$ gives $\lambda_{z'}=\lambda_x\eqqcolon\lambda$. The linear polynomials $F^w_z-\lambda F^{w'}_z$ consequently vanish
for $z$ equal to $x$ or one of its covers in $[x,y]$.
\Cref{lem:linear-section,lem:coatom-ratio} then show that they vanish
throughout the interval. At $z=y$ this equates nonzero multiples of the distinct
root lines of $y\to w$ and $y\to w'$, a contradiction.
\end{proof}

\section{Purity and connectedness of slices}\label{sec:geometry}

Throughout this section, $W$ will denote the crystallographic Coxeter group associated to an integral Cartan matrix $A$, with a fixed Kac--Moody realization of its root system. We write $\mathcal G$ for a complex Kac--Moody group associated with $A$
whose maximal torus $\mathcal T$ has Lie algebra
$\mathfrak h_\R\otimes\C$, and let $\mathcal B\supseteq\mathcal T$
be its positive Borel subgroup. Any choice of integral form of
$\mathfrak h_\R$ (with $\alpha_i^\vee$ integral and $\alpha_i$
integer-valued on it) gives such a group; these differ by central
isogenies inside $\mathcal T$, so $\mathcal G/\mathcal B$ and its
Schubert varieties $X_w=\overline{\mathcal{B}w\mathcal{B}/\mathcal{B}}$ do not depend on the choice. We refer the reader to \cite[Ch.~7]{Kumar} for background on $\mathcal{G/B}$.

Our goal is to prove the following proposition about the graphs $\mathcal M_I(u)$ of \Cref{def:slice}.  

\begin{proposition}\label{prop:slice-geometry}
Let $A$ be an integral Cartan matrix with Weyl group $W$. Let $U\subseteq V$ be spanned by roots, let
$H=H_U$, let $x\le w \in W$, and let $I=I_H(w;x)$ be the slice. Suppose $x=ux_0$ has
$q_w(x)=0$. Then $\ell_H(v)=d_I(u)$ for all $v\in\mathcal M_I(u)$ and $\mathcal M_I(u)$ is connected.
\end{proposition}
\begin{proof}
We have open, affine, connected, $\mathcal T$-stable neighborhood $Z=\dot{x}\mathcal{B}^-\mathcal{B}/\mathcal{B} \cap X_w$ 
of $x\mathcal{B}$ in $X_w$ such that some cocharacter
of $\mathcal T$ contracts $Z$ onto $x\mathcal{B}$.
Since $Z$ is nonempty and open in the irreducible
variety $X_w$, it is irreducible of dimension $L=\ell(w)$ and is
Cohen--Macaulay by \cite[Thm.~3(a)]{Mathieu}.

The $\mathcal T$-invariant curves through $x\mathcal{B}$ correspond to the edges incident to $x$ in $\Gamma_w$, with tangent space characters corresponding to their root labels. Attractivity implies that all
one-dimensional $\mathcal T$-orbits in $Z$ have $x\mathcal{B}$ in their
closure. Thus the number of curve coordinates in Brion's
finite-map construction \cite[Thm.~17]{BrionEquivariant},
in the form of \cite[Lem.~3.4]{FW}, is
$d_w(x)$. That construction gives a finite $\mathcal T$-equivariant morphism
\[
 f:Z\longrightarrow\A^{d_w(x)},
\]
with one coordinate for each invariant curve, whose character
is plus or minus the corresponding root. Since $q_w(x)=0$, we have $d_w(x)=L=\dim Z$. The coordinates of $f$ vanish at $x\mathcal{B}$ and, since $f$ is
finite, generate an ideal primary to the maximal ideal of the
$L$-dimensional Cohen--Macaulay local ring
$\mathcal O_{Z,x\mathcal{B}}$. 

The integral cocharacters annihilating $U$ are those of a subtorus
$\mathcal S\subseteq\mathcal T$, and a character of $\mathcal T$
is trivial on $\mathcal S$ exactly when it lies in $U$. By
\eqref{eq:root-reflections}, the positive roots trivial on
$\mathcal S$ are exactly those in $\Phi_H^+$. A coordinate of $\A^L$ is fixed by $\mathcal S$ exactly when it corresponds to an edge in the slice $I$. Hence the fixed coordinate subspace has dimension
$D=d_I(u)$. Its inverse image $Y=f^{-1}\bigl((\A^L)^{\mathcal S}\bigr)$ is the zero scheme of the other $L-D$ coordinates. The $D$ fixed
coordinates cut $\mathcal O_{Y,x\mathcal{B}}$ down to dimension zero,
so it has dimension $D$ and is Cohen--Macaulay. Moreover, $Y_{\mathrm{red}}=(Z^{\mathcal S})_{\mathrm{red}}$ by \cite[proof of Lem.~3.4]{FW}.

The local scheme $\operatorname{Spec}\mathcal O_{Y,x\mathcal{B}}$ is
connected and Cohen--Macaulay, and all its irreducible components
have dimension $D$. By Hartshorne's connectedness theorem
\cite[Exp.~III, Thm.~3.6 and Cor.~3.7]{SGA2}, any two components can be joined by a sequence of components whose consecutive intersections have dimension $D-1$. When $D=0$, the local scheme
has only one point and one irreducible component. Passing to
the reduction changes neither these components nor the
dimensions of their intersections.

We identify this local component structure with the intrinsic
Schubert geometry of the slice. For
$h\in H$ let $\operatorname{Inv}_H(h)=\{\beta\in\Phi_H^+ \mid h^{-1}\beta\in\Phi_H^-\},$
a set of size $\ell_H(h)$. Write $\Omega_z=\mathcal{B}z\mathcal{B}/\mathcal{B}$. Use ordered inversion-root coordinates
on $\Omega_{hx_0}$, in any order \cite[Ch.~7]{Kumar}, and write
$\mathcal U_\beta$ for the root subgroup of $\beta$. The
$\mathcal S$-fixed equations set the coordinates with root outside
$U$ to zero. The remaining coordinate roots are
\[
 \{\beta\in\Phi_H^+ \mid (hx_0)^{-1}\beta\in\Phi^-\}
 =\operatorname{Inv}_H(h),
\]
since $x_0^{-1}(\Phi_H^+)\subseteq\Phi^+$, the defining property of
the coset minimum $x_0$. Hence
\[
 C_h\coloneqq(\Omega_{hx_0}^{\mathcal S})_{\mathrm{red}}
 =\Bigl(\prod_{\beta\in\operatorname{Inv}_H(h)}\mathcal U_\beta\Bigr)hx_0\mathcal{B}
\]
is an affine space of dimension $\ell_H(h)$.

We claim that for all $v \in H$ we have
\begin{equation}\label{eq:slice-closure}
 \overline{C_v}=\bigsqcup_{h\le_Hv}C_h.
\end{equation}
The groups $\mathcal T$ and $\mathcal U_\gamma$ for
$\gamma\in\Phi_H^+$ lie in $\mathcal{B}$ and commute with
$\mathcal S$, so they preserve each $C_h$. Let
$s=s_\beta\in S_H$ with $\beta\in\Phi_H^+$, let $\mathcal P_\beta$
be generated by $\mathcal T$ and $\mathcal U_{\pm\beta}$, and set
$\mathcal B_\beta=\mathcal T\mathcal U_\beta$. Then
$\mathcal P_\beta/\mathcal B_\beta\cong\mathbb P^1$ and
$\mathcal P_\beta=\mathcal B_\beta\sqcup\mathcal U_\beta \dot{s}\mathcal B_\beta$. Suppose
$sh>_Hh$. As $s$ permutes $\Phi_H^+\setminus\{\beta\}$, we have
$\operatorname{Inv}_H(sh)=\{\beta\}\sqcup s\operatorname{Inv}_H(h)$,
so $C_{sh}=\mathcal U_\beta \dot{s}C_h$ and therefore
\[
 \mathcal P_\beta C_h=\mathcal P_\beta C_{sh}=C_h\sqcup C_{sh}.
\]
Now induct on $\ell_H(v)$, starting from $C_e=\{x_0\mathcal{B}\}$.
Write $v=sv'$ with $\ell_H(v')=\ell_H(v)-1$. By induction,
$\overline{C_{v'}}=\bigsqcup_{h\le_Hv'}C_h$ is irreducible and
$\mathcal B_\beta$-stable, so the map
$\mathcal P_\beta\times^{\mathcal B_\beta}\overline{C_{v'}}\to\mathcal{G}/\mathcal{B}$
has complete irreducible source. Its image
$\mathcal P_\beta\overline{C_{v'}}$ is therefore closed and
irreducible, and it equals $\bigsqcup_{h\le_Hv}C_h$ by the subword
property of Bruhat order. Since $C_v$ is its unique cell of top dimension, this
proves \eqref{eq:slice-closure}.

Every point of $Z^{\mathcal S}$ lies in $C_h$ for some $h\in I$.
Indeed, let $p\in Z^{\mathcal S}\cap\Omega_z$, so that $z\le w$. The
closure of $\mathcal Tp$ is a projective toric variety fixed
pointwise by $\mathcal S$. It contains $x\mathcal{B}$ and
$z\mathcal{B}$, the limits of $p$ under cocharacters contracting
$Z$ and $\Omega_z$. Its $\mathcal T$-fixed points are joined by
$\mathcal T$-invariant curves, whose roots lie in $U$ by
\cite[Prop.~7.2, p.~525]{FW} because $\mathcal S$ fixes these curves
pointwise. Thus $z\in Hx$.

Since $I$ is a finite intrinsic order ideal,
\eqref{eq:slice-closure} now gives
\[
 (Z^{\mathcal S})_{\mathrm{red}}
 =Z\cap\bigcup_{v\in I}C_v
 =Z\cap\bigcup_{v\in \mathcal{M}_I(u)}\overline{C_v}.
\]
The varieties in the final union are irreducible and none
contains another. Such a variety contains $x\mathcal{B}=ux_0\mathcal{B}$
exactly when $u\le_Hv$. Therefore the irreducible components of
$Z^{\mathcal S}$ through $x\mathcal{B}$ are the $Z\cap\overline{C_v}$
with $v\in\mathcal M_I(u)$. Each is nonempty open, and hence dense,
in $\overline{C_v}$.
Its dimension and its local dimension at the closed point $x\mathcal{B}$
are both $\ell_H(v)$. Purity of $Y$ now gives $\ell_H(v)=D$ for $v\in\mathcal M_I(u)$.

Finally, by \eqref{eq:slice-closure}, an intersection
$\overline{C_v}\cap\overline{C_{v'}}$ is the finite union of the
$\overline{C_h}$ over common lower bounds $h$ of $v$ and $v'$. Its local dimension at $x\mathcal{B}$ is
therefore the largest $\ell_H(h)$ among these bounds satisfying
$u\le_Hh$. If the local intersection has dimension $D-1$,
there is such an $h$ with $\ell_H(h)=D-1$. As both $v$ and
$v'$ have length $D$, this means
\[
 u\le_Hh,\qquad h\lessdot_Hv,\qquad h\lessdot_Hv'.
\]
Every consecutive pair in the component sequence supplied by
connectedness is thus an edge of $\mathcal M_I(u)$. This proves $\mathcal{M}_I(u)$ is a connected graph, and completes the proof.
\end{proof}

The map $f$ from the above proof can also be used to bound Bruhat graph degrees in $I$.

\begin{lemma}\label{lem:strict-length}
Let $A$ be an integral Cartan matrix with Weyl group $W$. Let $x\le w \in W$
with $q_w(x)=0$, and let $U\subseteq V$ be spanned by roots
with $\dim U<\lvert\supp(w)\rvert$. Set $H=H_U$, let $x_0$ be the minimum
of $Hx$, and write $x=ux_0$. For the slice
$I=I_H(w;x)$, one has
\[
 d_I(u)<\ell(w).
\]
\end{lemma}
\begin{proof}
By \eqref{eq:root-reflections}, $d_I(u)$ counts the roots
$\gamma\in\Phi^+\cap U$ with $s_\gamma x\le w$. Hence
$d_I(u)\le d_w(x)=\ell(w)$, with equality only if every root
$\gamma\in\Phi^+$ with $s_\gamma x\le w$ lies in $U$. Suppose this
happens, and let $Z$, $f$, and $\mathcal S$ be as in the proof of
\Cref{prop:slice-geometry}. The character associated to each coordinate of $f$
is a root in $U$,
so $\mathcal S$ acts trivially on the target of $f$. As $f$ is
$\mathcal S$-equivariant with finite fibers and $\mathcal S$ is
connected, $\mathcal S$ fixes $Z$ pointwise, and hence also its
closure $X_w$. For $s\in\supp(w)$, the curve
$X_s\subseteq X_w$ joins $e\mathcal B$ and $s\mathcal B$, and
$\mathcal T$ acts on it with character $\pm\alpha_s$. As
$\mathcal S$ fixes $X_s$ pointwise, $\alpha_s\in U$. Thus $U$
contains the $\lvert\supp(w)\rvert$ linearly independent simple roots
$\alpha_s$ with $s\in\supp(w)$, contradicting
$\dim U<\lvert\supp(w)\rvert$.
\end{proof}

\section{Proof of the Lookup Conjecture}\label{sec:reduction}

We are now ready to prove \Cref{thm:main} and thereby the Lookup Conjecture.

\begin{proof}[Proof of \Cref{thm:main}]
Let $W$ be the Weyl group of an integral Cartan matrix, and suppose that the Lookup Conjecture holds in every reflection subgroup of $W$ of intrinsic rank three. Suppose for the sake of contradiction that the Lookup Conjecture fails in some reflection subgroup of $W$, and choose a counterexample $(\check{x},w)$ in a reflection subgroup $W'$ such that the intrinsic length $\ell_{W'}(w)$ is minimal over all counterexamples in all reflection subgroups; call this minimal length $L$. Let $\supp_{W'}(w) \subset S'$ denote the support of $w$ in the intrinsic Coxeter structure $(W',S')$, and let $K$ be the standard parabolic subgroup of $W'$ generated by $\supp_{W'}(w)$. Since the Bruhat graph below $w$ in $K$ is equal to that in $W'$, it is immediate that $(\check{x},w)$ is also a counterexample to the Lookup Conjecture in $K$ (and $\ell_K(w)=L$ is still minimal). By hypothesis, since $K \leq W' \leq W$ is a reflection subgroup of $W$ of intrinsic rank $|\supp_{W'}(w)|$, we have that $|\supp_{W'}(w)|>3$ (a counterexample in a dihedral group is impossible by \Cref{prop:easy-dihedral}).

Fix a Kac--Moody realization of the root system of $K$. From now on we work inside $K$, and symbols without subscripts refer to the intrinsic quantities in $K$.

Since $(\check{x},w)$ is a counterexample to (\ref{eq:lookup-hypothesis}) implying (\ref{eq:CP-criterion}), it must be that there are $x \xrightarrow{\alpha} y$ in $[\check{x},w]$ such that $x$ satisfies (\ref{eq:lookup-hypothesis}), with respect to $w$, but $y$ does not. We have $q_w(y)=0$, so there must be $y\xrightarrow{\beta} z\le w$ with $q_w(z)>0$.  By \Cref{prop:plane-comparison}, there is $M \in \mathcal{D}_\beta$ with
$\delta_M(y,z)=1$. Choose $\gamma\in\Phi^+$ with
$s_\gamma\in M\setminus\{s_\beta\}$, and set
\[
 U=\Span_\R\{\alpha,\beta,\gamma\},\qquad H=H_U.
\]
The roots of $M$ lie in
$\Span_\R\{\beta,\gamma\}\subseteq U$, so $M\subseteq H$.

Take the slice $I=I_H(w;x)$ on $Hx$. Write
$x=ux_0$, $y=u_1x_0$, and $z=u_2x_0$, so that $u\to u_1\to u_2$, and set $D=d_I(u)$. Since $\dim{U} \leq 3 < |\supp(w)|$, \Cref{lem:strict-length} gives
\begin{equation}\label{eq:induction-decrease}
 D<L.
\end{equation}
\Cref{prop:slice-geometry} gives
$\ell_H(v)=D$ for every $v\in\mathcal M_I(u)$, and
connectedness of $\mathcal M_I(u)$. Moreover, since defects are nonnegative \cite{DyerDeodhar} and since $[e,v]_H \subseteq I$, we have
\begin{equation}
\label{eq:all-neighbors-below-v}
D=\ell_H(v)\le d_{[e,v]_H}(u)\le d_I(u)=D.
\end{equation}
Thus every reflection neighbor $u'\leftarrow u$ in $I$ lies in
every $[e,v]_H$, and since $x \to u'x_0$ we have $q_w(u'x_0)=0$. Thus by \Cref{cor:internal} we have $d_I(u')=d_I(u)=D$.
As in (\ref{eq:all-neighbors-below-v}), we see
\begin{equation}\label{eq:component-neighbors}
 d_{[e,v]_H}(u')=D,\qquad
 N_{[e,v]_H}(u')=N_I(u').
\end{equation}
Thus $u$ satisfies (\ref{eq:lookup-hypothesis}) in $[e,v]_H$.

Since $H$ is a reflection subgroup of $W$ and $\ell_H(v)=D<L$, minimality of $L$ implies that $[e,v]_H$ is regular above $u$ for all $v \in \mathcal M_I(u)$. Moreover, if $\mathcal M_I(u)$ had more than one element, the connectedness of \Cref{prop:slice-geometry} would give distinct $v,v'$ with a common lower cover $h\ge_Hu$. By (\ref{eq:component-neighbors}),
$N_{[e,v]_H}=N_I=N_{[e,v']_H}$ at $u$ and at every upward
reflection neighbor of $u$ in $I$, contradicting \Cref{prop:separation}. Hence $\mathcal M_I(u)=\{v\}$ is a singleton. For any $h\in I$ with $u\le_Hh$, we thus have $N_I(h)=N_{[e,v]_H}(h)$. And regularity of $[e,v]_H$ above $u$ now gives $d_I(h)=D$. In particular, 
$d_I(u_1)=d_I(u_2)=D$. On the other hand, $q_w(y)=0$, $M\subseteq H$, and $\delta_M(y,z)=1$. \Cref{cor:internal} therefore forces $d_I(u_2)>d_I(u_1)$, a contradiction. 
\end{proof}

\subsection{Finite and affine groups}\label{sec:corollaries}

We now apply \Cref{thm:main} to our cases of interest.

\begin{proof}[Proof of \Cref{cor:finite}]
First let $W$ be finite and crystallographic. Its reflection
subgroups are finite crystallographic Coxeter groups. Their
irreducible rank-three possibilities are $A_3$ and $B_3$;
the other possibilities are products of rank-one and dihedral
groups. Boe--Graham's finite Weyl-group checks \cite[p.~319]{BG} supply
the hypothesis of \Cref{thm:main}.
\end{proof}

\begin{proof}[Proof of \Cref{cor:affine}]
Let $W$ be the Weyl group of type $\widetilde{A}_n$. By the Dyer--Lehrer classification
\cite[Thm.~3]{DL}, all reflection subgroups of $W$ are products of groups of types $A$ and $\widetilde{A}$. So an irreducible reflection subgroup of intrinsic rank three is therefore of type $A_3$ or $\widetilde A_2$. The Lookup Conjecture for $\widetilde A_2$ was proven by Boe--Graham \cite{BGAffine} and Graham--Li \cite{GL}, so we may apply \Cref{thm:main}.
\end{proof}

\section{Degree monotonicity in parabolic quotients}
\label{sec:grassmannian}
Recall that for $J\subseteq S$ and $w\in W^J$, the \emph{parabolic Bruhat graph} $\Gamma_w^J$ has vertex set $[e,w]\cap W^J$, and edges $x \to y$ if $tx W_J=yW_J$ for some $t\in T$ and $\ell(y)>\ell(x)$. The degree of a vertex $x$ in this graph is denoted $d_w^J(x)$. In this section, we analyze the monotonicity of $d_w^J$ for certain important parabolic quotients. In all the cases considered in this section, $J$ generates a finite Coxeter group $W_J$, and we write $w_0(J)$ for its longest element.
\begin{lemma}\label{lem:degree-quotient-to-full}
For $x^J\leq w^J$ in $W^J$, $d_{w^J}^J(x^J)+\ell(w_0(J))=d_{w^Jw_0(J)}(x^Jx_J)$ for any $x_J\in W_J$.
\end{lemma}
\begin{proof}
Write $T_J$ for $W_J \cap T$, and let $x=x^Jx_J$. We first show that $xtW_J=xt'W_J$ for $t\neq t'\in T$ if and only if both $t,t'\in T_J$. The ``$\Leftarrow$" direction is clear. For the ``$\Rightarrow$" direction, assume that $xtW_J=xt'W_J$, so $t=t'z$, $z=t't$ for some $z\in W_J$. Consider the reflection subgroup $D$ generated by $t$ and $t'$, and let $a$ and $b$ be the intrinsic simple generators of $D$. Since $z\in D$ and $z$ has even length, $z\neq a,b$ and thus $z>a,b$ in the Bruhat order. Since $W_J$ is an order ideal, we conclude $a,b,t,t' \in W_J$. 

Now, $d_{w^Jw_0(J)}(x^Jx_J)$ counts the number of reflections $t$ such that $x^Jx_J t\leq w^Jw_0(J)$. Any of the $\ell(w_0(J))$ reflections $t\in T_J$ works, as $x^J\leq w^J$ and $x_Jt\leq w_0(J)$. And we have just seen that each $t\in T\setminus T_J$ produces a different coset representative $xt$, and the number of such edges is $d_{w^J}^J(x^J)$ by definition.
\end{proof}

\Cref{lem:degree-quotient-to-full} allows us to study degree monotonicity inside a parabolic quotient by lifting to the full Coxeter group. 

\subsection{Degree monotonicity for affine Grassmannians}
In this subsection, $W$ has type $\widetilde A_{n-1}$ with generators $S=\{s_0,s_1,\ldots,s_{n-1}\}$, with the convention that $s_0$ denotes the affine node; set $J=\{s_1,\ldots,s_{n-1}\}$. Correspondingly, let  $\Phi_{\mathrm{fin}}$ be the finite root system of $W_J$ with simple roots $\alpha_1,\ldots\alpha_{n-1}$ and let $\delta=\alpha_0+\alpha_1+\cdots+\alpha_{n-1}$ be the null root. Then each root of $W$ is written as $\beta+k\delta$ with $\beta\in\Phi_{\mathrm{fin}}$ and $k\in\Z$, with the positive roots given by \[\Phi^+=\{\beta+k\delta\:|\: \beta\in\Phi_{\mathrm{fin}},k\in\Z_{>0}\}\cup\Phi^+_{\mathrm{fin}}.\]
We will also use the realization of $W$ as the group of \emph{affine permutations}: \[W=\left\{\text{bijections }w:\Z\rightarrow\Z\:\bigg|\: w(i{+}n)=w(i)+n\text{ for all }i\in\Z,\text{ and}\ \sum_{i=1}^n w(i)={n+1\choose 2}\right\}.\]
In this realization, we have reflection set $T=\{t_{a,b}\:|\:a\not\equiv b\bmod n,a<b,a\in[n]\}$, where $t_{a,b}(a+kn)=b+kn$ and $t_{a,b}(b+kn)=a+kn$ for $k\in\Z$, and $t_{a,b}(i)=i$ if $i\not\equiv a,b\bmod n$. 

The following is the \emph{affine tableau criterion} for Bruhat order on $W$ (see \cite[Thm.~8.3.7]{BB}). 
\begin{proposition}\label{prop:tableau-affine-Bruhat}
For $w\in W(\widetilde A_{n-1})$ and $p,q\in\Z$, define the finite rank count \[r_w(p,q)=|\{i\leq p\:|\:w(i)\geq q\}|.\]
Then $x\leq y$ if and only if $r_x(p,q)\leq r_y(p,q)$ for all $p,q\in\Z$.
\end{proposition}

For a $2$-dimensional subspace $P \subset V$ spanned by two roots, the subgroup $H_P$ has type $A_1\times A_1$, $A_2$, or $\widetilde A_1$. Specifically, let $\beta,\gamma$ be a basis of $P$ with corresponding reflections $t_{a,b}$ and $t_{c,d}$. Consider the cardinality $m=|\{a\bmod n, b\bmod n,c\bmod n,d\bmod n\}|$. If $m=4$, $t_{a,b}$ and $t_{c,d}$ commute, so $H_P$ is of type $A_1\times A_1$. If $m=3$, then $H_P$ is of type $A_2$. If $m=2$, then $\delta\in P$ so we can choose $\beta$ and $\gamma$ so that $\beta+\gamma=\delta$, and $H_P$ is of type $\widetilde A_1$.

\begin{proof}[Proof of \Cref{thm:grassmannian}]
Fix $x^J\lessdot y^J\leq w^J$ in $W^J$. Let $x=x^Jw_0(J)$, $y=xs_{\alpha}=y^Jy_J$ where $s_{\alpha}\in T\setminus T_J$, and $w=w^Jw_0(J)$. By \Cref{lem:degree-quotient-to-full}, we need to show that $d_w(x)\geq d_w(y)$. In this proof we associate edges of the Bruhat graph to the corresponding \emph{right} reflection, in contrast to the left-sided conventions of earlier sections.

Let $d_{M,w}(y)$ and $d_{M,w}(x)$ be as in \Cref{sec:degree-comparison}, where $M\in\mathcal{D}_{\alpha}$ is a maximal dihedral reflection subgroup containing $s_{\alpha}$. Let $z\in W$ be the minimal element of the coset $xM=yM$. By \eqref{eq:plane-sum}, it suffices to show that $d_{M,w}(y)\leq d_{M,w}(x)$, which is clear if $M$ is of type $A_1\times A_1$.

If instead $M$ is of type $\widetilde A_1$, we claim that $zM\cap [e,w]$ contains a unique maximum. As explained above, when $M$ is of type $\widetilde A_1$, we can pick intrinsic generators $s_{\beta}$ and $s_{\gamma}$ such that $\beta,\gamma\in\Phi^+$, $\beta+\gamma=\delta=\alpha_0+\alpha_1+\cdots+\alpha_{n-1}$. Without loss of generality, assume that $\alpha_0$ is supported on $\gamma$, and thus $\beta\in\Phi_{\mathrm{fin}}$, $s_{\beta}\in W_J$. If $zM\cap [e,w]$ were to contain two maximal elements $z_1=za_1$ and $z_2=za_2$ of intrinsic length $\ell_M(a_1)=\ell_M(a_2)=\ell$, then we can assume $a_1s_{\beta}$ and $a_2s_{\gamma}$ are the two elements of $M$ of intrinsic length $\ell+1$. Since $za_1\leq w$ and $s_{\beta}\in W_J$, we have $za_1s_{\beta}\leq w$ as the parabolic subgroup part of $w=w^Jw_0(J)$ is maximal. Both $za_1,za_2< za_1s_{\beta}\leq w$, a contradiction. Because $zM\cap [e,w]$ has a unique maximum, the degree function $d_{M,w}$ in this dihedral coset is constant, so $d_{M,w}(x)=d_{M,w}(y)$.

Finally, if $M$ is of type $A_2$ then $zM$ contains $6$ elements. The only possibility allowing for $d_{M,w}(y)>d_{M,w}(x)$ is that $zM\cap[e,w]=\{z,z_1,z_2,z_3,z_4\}$ where $z< z_i< z_j$ for $i\in\{1,2\}$ and $j\in\{3,4\}$, in which case $x=z$ and $y\in\{z_1,z_2\}$ would give $d_{M,w}(y)=3$ and $d_{M,w}(x)=2$. We show that this cannot happen, by deducing that the maximum $z_5$ of $zM$ satisfies $z_5\leq w$ whenever $z_3,z_4\leq w$, using \Cref{prop:tableau-affine-Bruhat}.

Write $z_1=zt_{a,b}$ and $z_2=zt_{b,c}$, where $a,b,$ and $c$ are all distinct modulo $n$. Then $a<b<c$ and $z(a)<z(b)<z(c)$. Assume without loss of generality that $z_3(c)<z_3(a)<z_3(b)$ and $z_4(b)<z_4(c)<z_4(a)$. We claim that if $p-q\geq b-z(b)$ then $r_{z_5}(p,q)=r_{z_3}(p,q)$, and if $p-q<b-z(b)$ then $r_{z_5}(p,q)=r_{z_4}(p,q)$. Fix $p,q\in\Z$ with $p-q\geq b-z(b)$, and count the contributions of the coordinates $a+kn,b+kn,$ and $c+kn$ to $r_{z_3}(p,q)$ and $r_{z_5}(p,q)$ for each $k\in\Z$. Specifically, for $u\in\{z_3,z_5\}$, we count the number of $i\in\{a+kn,b+kn,c+kn\}$ such that $i\leq p$ and $u(i)\geq q$. It is routine to check that if $p<a+kn$ or $q>z(c)+kn$, both contributions are $0$; if $a+kn\leq p<b+kn$ or $z(b)+kn<q\leq z(c)+kn$, both contributions are $1$; if $ b+kn\leq p<c+kn$ and $z(a)+kn<q\leq z(b)+kn$, both contributions are $2$; and if $p\geq c+kn$ and $q\leq z(a)+kn$, both contributions are $3$. See \Cref{fig:degree-monotone-A2} for a visualization, where our point $(p,q)$ is chosen weakly below the line with slope $1$ that passes through $(b,z(b))$, and the rank count $r(p,q)$ is the number of dots weakly northwest of the coordinate $(p,q)$. Checking $r_{z_4}(p,q)=r_{z_5}(p,q)$ when $p-q<b-z(b)$ is exactly analogous. As a result, for all $p,q\in\Z$, \[r_{z_5}(p,q)=\max(r_{z_3}(p,q),r_{z_4}(p,q))\leq r_w(p,q)\]
so $z_5\leq w$, a contradiction. We conclude that $d_{M,w}(y)\leq d_{M,w}(x)$. 
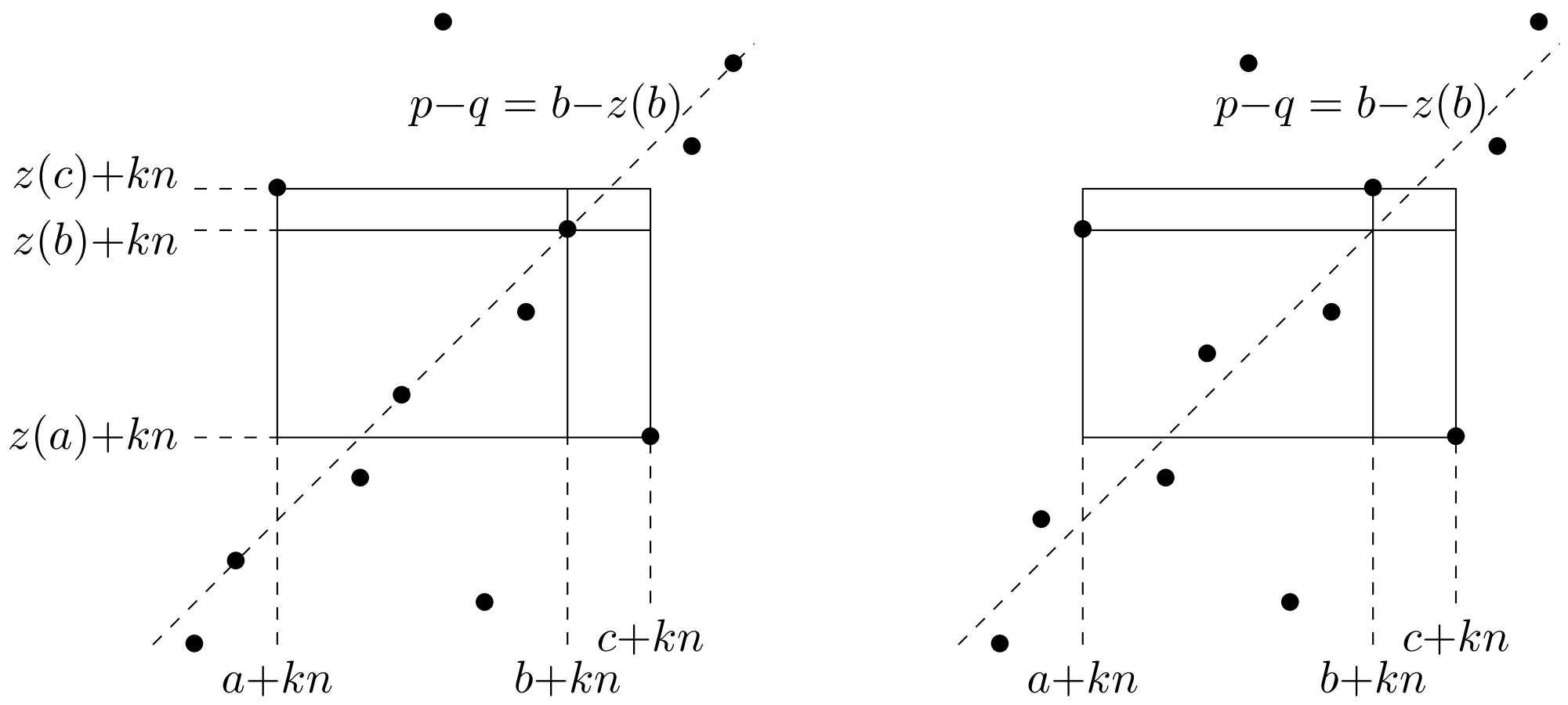
\begin{figure}[h!]
\centering
\begin{tikzpicture}[scale=0.35]
\node at (4,11) {$\bullet$};
\node at (8,15) {$\bullet$};
\node at (2,0) {$\bullet$};
\node at (6,4) {$\bullet$};
\node at (10,8) {$\bullet$};
\node at (14,12) {$\bullet$};
\node at (3,2) {$\bullet$};
\node at (7,6) {$\bullet$};
\node at (11,10) {$\bullet$};
\node at (15,14) {$\bullet$};
\node at (9,1) {$\bullet$};
\node at (13,5) {$\bullet$};
\draw(4,5)--(13,5)--(13,11)--(4,11)--(4,5);
\draw(11,5)--(11,11);
\draw(4,10)--(13,10);
\node[below] at (4,0) {$a{+}kn$};
\node[below] at (11,0) {$b{+}kn$};
\node[below] at (13,1) {$c{+}kn$};
\draw[dashed](4,0)--(4,5);
\draw[dashed](11,0)--(11,5);
\draw[dashed](13,1)--(13,5);
\node[left] at (2,5) {$z(a){+}kn$};
\node[left] at (2,9.7) {$z(b){+}kn$};
\node[left] at (2,11.3) {$z(c){+}kn$};
\draw[dashed](2,5)--(4,5);
\draw[dashed](2,10)--(4,10);
\draw[dashed](2,11)--(4,11);
\draw[dashed](1,0)--(15.5,14.5);
\node at (10.5,13) {$p{-}q=b{-}z(b)$};
\end{tikzpicture}
\qquad\qquad
\begin{tikzpicture}[scale=0.35]
\node at (4,10) {$\bullet$};
\node at (8,14) {$\bullet$};
\node at (2,0) {$\bullet$};
\node at (6,4) {$\bullet$};
\node at (10,8) {$\bullet$};
\node at (14,12) {$\bullet$};
\node at (3,3) {$\bullet$};
\node at (7,7) {$\bullet$};
\node at (11,11) {$\bullet$};
\node at (15,15) {$\bullet$};
\node at (9,1) {$\bullet$};
\node at (13,5) {$\bullet$};
\draw(4,5)--(13,5)--(13,11)--(4,11)--(4,5);
\draw(11,5)--(11,11);
\draw(4,10)--(13,10);
\node[below] at (4,0) {$a{+}kn$};
\node[below] at (11,0) {$b{+}kn$};
\node[below] at (13,1) {$c{+}kn$};
\draw[dashed](4,0)--(4,5);
\draw[dashed](11,0)--(11,5);
\draw[dashed](13,1)--(13,5);
\draw[dashed](1,0)--(15.5,14.5);
\node at (10.5,13) {$p{-}q=b{-}z(b)$};
\end{tikzpicture}
\caption{For visualization in the proof of \Cref{thm:grassmannian} in the case of type $A_2$ dihedral coset, $z_5$ is on the left and $z_3$ on the right.}
\label{fig:degree-monotone-A2}
\end{figure}
\end{proof}

\begin{proof}[Proof of \Cref{cor:grassmannian-lookup}]
Suppose that $x \leq z \leq w$ with $x,w \in M_J$ and that $q_w(x)=0$. Then $x^J \leq z^J \leq w^J$, so by \Cref{thm:grassmannian} we have $d^J_{w^J}(x^J) \geq d^J_{w^J}(z^J)$. Since $w\in M_J$, we have $w=w^Jw_0(J)$, so by \Cref{lem:degree-quotient-to-full} we have $d_w(x) \geq d_w(z)$ hence $0=q_w(x) \geq q_w(z) \geq 0$. Thus $q_w(z)=0$.
\end{proof}

\subsection{Weight diagrams and geometric Satake}
We now discuss a consequence of \Cref{thm:grassmannian}.

Write $W_J$ for the finite Weyl group of $\mathcal{G}=\operatorname{SL}_{n+1}(\C)$, a parabolic subgroup of the type $A$ affine Weyl group $W$, and let $w_\nu$ denote the unique Bruhat-maximal element of $W_J\tau^\nu W_J$. Let $Y_\lambda = X^J_{w_\lambda^J} \subset \operatorname{Gr}_{\mathcal{G}}$ denote the spherical Schubert variety indexed by a dominant coweight $\lambda$, and $V_\lambda$ the corresponding irreducible representation of $\mathcal{G}^\vee=\operatorname{PGL}_{n+1}(\C)$.
Write $\Psi^+$ for the positive roots of $\mathcal{G}^\vee$ and
$\operatorname{Wt}(V_\lambda)$ for the set of weights of $V_\lambda$. For
$\eta\in\operatorname{Wt}(V_\lambda)$, define the quantity
\[
 R_\lambda(\eta)\coloneqq
 \sum_{\beta\in\Psi^+}
 \lvert\{k\in\Z\setminus\{0\} \mid
            \eta+k\beta\in\operatorname{Wt}(V_\lambda)\}\rvert.
\]

\begin{corollary}\label{cor:weight-count}
Let $\mu\le \nu \leq \lambda$ be dominant coweights of $\mathcal{G}$.
Then
\begin{equation}
\label{eq:R-minus-dim}
R_\lambda(\mu)-\dim Y_\lambda=q_{w_\lambda}(w_\mu).
\end{equation}
Therefore $R_\lambda(\mu)\ge R_\lambda(\nu)\ge\dim Y_\lambda$ and
$R_\lambda(\mu)=\dim Y_\lambda$ if and only if
$\dim V_\lambda[\mu]=1$.
\end{corollary}
\begin{proof}
The $\mathcal{T}$-fixed points of $Y_\lambda$ are precisely
$\tau^\eta \mathcal{G}(\mathcal O)$ with
$\eta\in\operatorname{Wt}(V_\lambda)$
(see, e.g. \cite[Rem.~3.3]{MV}) and $\ell(w_\lambda^J)=\dim Y_\lambda$. For a positive root $\alpha$ of $\mathcal{G}$ and $m\in\Z$,
the affine reflection in the hyperplane
$\langle\alpha,\cdot\rangle=m$ sends $\eta\longmapsto
 \eta+(m-\langle\alpha,\eta\rangle)\alpha^\vee$. Consequently $d_{w_\lambda^J}^J((\tau^\eta)^J)=R_\lambda(\eta)$. The count $R_\lambda$ is invariant under $W_J$.
Moreover, $(w_\xi)^J=(\tau^{y\xi})^J$ for some
$y\in W_J$. Thus \Cref{lem:degree-quotient-to-full} gives $q_{w_\lambda}(w_\mu) 
 =R_\lambda(\mu)-\ell(w_\lambda^J)$, proving (\ref{eq:R-minus-dim}).

Dominance implies $w_\mu\le w_\nu\le w_\lambda$
\cite[\S2]{LusztigWeights}.
By \Cref{thm:grassmannian} and \Cref{lem:degree-quotient-to-full} we have $R_\lambda(\mu) \geq R_{\lambda}(\nu)$, and the nonnegativity of defects \cite{DyerDeodhar} implies $R_{\lambda}(\nu) \geq \dim(Y_\lambda)$ by (\ref{eq:R-minus-dim}).

Finally, $\dim V_\lambda[\mu]=P_{w_\mu,w_\lambda}(1)$ by
\cite[Thm.~6.1]{LusztigWeights}.
This value is one precisely when $Y_\lambda$ is rationally
smooth at $\tau^\mu \mathcal{G}(\mathcal O)$
\cite[Rem.~2.3]{HainesRicharz}.
\Cref{cor:grassmannian-lookup} therefore implies $\dim V_\lambda[\mu]=1$ if and only if $q_{w_\lambda}(w_\mu)=0$.
\end{proof}

\subsection{Degree monotonicity in finite types}\label{sub:degree-monotone-finite-type}
We say that the parabolic quotient $W^J$ is \emph{degree-monotone} if for all $x\leq y\leq w\in W^J$ we have $d_w^J(x)\geq d_w^J(y)$.
\begin{lemma}
\label{lem:W-K-degree-monotone}
Let $(W,S)$ be a Coxeter system and suppose $W_K$ is finite. If $W^J$ is degree-monotone for some $J\subseteq K$, then $W^K$ is degree-monotone.  
\end{lemma}
\begin{proof}
Take $x\leq y\leq w\in W^K\subseteq W^J$. By \Cref{lem:degree-quotient-to-full}, \[d_w^K(x)=d_{w^Kw_0(K)}(x)-\ell(w_0(K))=d^J_{w^Kw_0(K)w_0(J)}(x)+\ell(w_0(J))-\ell(w_0(K)).\]
Since $W^J$ is degree-monotone, $d^J_{w^Kw_0(K)w_0(J)}(x)\geq d^J_{w^Kw_0(K)w_0(J)}(y)$ so $d_w^K(x)\geq d_w^K(y)$. 
\end{proof}

\begin{figure}[h!]
\centering
\begin{tikzpicture}[
    scale=0.7,
    vertex/.style={draw,shape=circle,fill=black,scale=0.5},
    box/.style={draw,shape=rectangle,fill=white,scale=1.2}
]

% A_n
\draw (0,0)--(3,0);
\draw (5,0)--(8,0);
\foreach \x in {0,1,2,3,5,6,7,8}
    \node[box] at (\x,0) {};
\foreach \x/\lab in {
    0/1,1/2,2/3,3/{},5/{},6/{},7/{n{-}1},8/n
}
    \node[vertex,label=below:{$\lab$}] at (\x,0) {};
\node at (4,0) {$\cdots$};
\node at (-1,0) {$A_n$};

% B_n
\draw (0,-2)--(3,-2);
\draw (5,-2)--(7,-2);
\draw (7,-2.04)--(8,-2.04);
\draw (7,-1.96)--(8,-1.96);
\draw (7.4,-1.85)--(7.6,-2)--(7.4,-2.15);
\node[box] at (0,-2) {};
\foreach \x/\lab in {
    0/1,1/2,2/3,3/{},5/{},6/{},7/{n{-}1},8/n
}
    \node[vertex,label=below:{$\lab$}] at (\x,-2) {};
\node at (4,-2) {$\cdots$};
\node at (-1,-2) {$B_n$};

% C_n
\draw (0,-4)--(3,-4);
\draw (5,-4)--(7,-4);
\draw (7,-4.04)--(8,-4.04);
\draw (7,-3.96)--(8,-3.96);
\draw (7.6,-3.85)--(7.4,-4)--(7.6,-4.15);
\node[box] at (8,-4) {};
\foreach \x/\lab in {
    0/1,1/2,2/3,3/{},5/{},6/{},7/{n{-}1},8/n
}
    \node[vertex,label=below:{$\lab$}] at (\x,-4) {};
\node at (4,-4) {$\cdots$};
\node at (-1,-4) {$C_n$};

% D_n
\draw (0,-6)--(3,-6);
\draw (5,-6)--(8,-6);
\draw (7,-6)--(7,-5.2);
\node[box] at (0,-6) {};
\node[box] at (8,-6) {};
\node[box] at (7,-5.2) {};
\foreach \x/\lab in {
    0/1,1/2,2/3,3/{},5/{},6/{},7/{n{-}2}
}
    \node[vertex,label=below:{$\lab$}] at (\x,-6) {};
\node[vertex,label=above:{$n{-}1$}] at (8,-6) {};
\node[vertex,label=left:{$n$}] at (7,-5.2) {};
\node at (4,-6) {$\cdots$};
\node at (-1,-6) {$D_n$};

% E_6
\draw (11,0)--(15,0);
\draw (13,0)--(13,0.8);
\node[box] at (11,0) {};
\node[box] at (15,0) {};
\foreach \x/\lab in {11/1,12/2,13/3,14/4,15/5}
    \node[vertex,label=below:{$\lab$}] at (\x,0) {};
\node[vertex,label=right:{$6$}] at (13,0.8) {};
\node at (10,0) {$E_6$};

% E_7
\draw (11,-2)--(16,-2);
\draw (13,-2)--(13,-1.2);
\node[box] at (16,-2) {};
\foreach \x/\lab in {11/1,12/2,13/3,14/4,15/5,16/6}
    \node[vertex,label=below:{$\lab$}] at (\x,-2) {};
\node[vertex,label=right:{$7$}] at (13,-1.2) {};
\node at (10,-2) {$E_7$};

% E_8
\draw (11,-4)--(17,-4);
\draw (13,-4)--(13,-3.2);
\foreach \x/\lab in {11/1,12/2,13/3,14/4,15/5,16/6,17/7}
    \node[vertex,label=below:{$\lab$}] at (\x,-4) {};
\node[vertex,label=right:{$8$}] at (13,-3.2) {};
\node at (10,-4) {$E_8$};

% F_4
\draw (11,-6)--(12,-6);
\draw (13,-6)--(14,-6);
\draw (12,-6.04)--(13,-6.04);
\draw (12,-5.96)--(13,-5.96);
\draw (12.4,-5.85)--(12.6,-6)--(12.4,-6.15);
\foreach \x/\lab in {11/1,12/2,13/3,14/4}
    \node[vertex,label=below:{$\lab$}] at (\x,-6) {};
\node at (10,-6) {$F_4$};

% G_2
\draw (17,-6)--(18,-6);
\draw (17,-6.05)--(18,-6.05);
\draw (17,-5.95)--(18,-5.95);
\draw (17.4,-5.85)--(17.6,-6)--(17.4,-6.15);
\node[vertex,label=below:{$1$}] at (17,-6) {};
\node[vertex,label=below:{$2$}] at (18,-6) {};
\node at (16,-6) {$G_2$};

\end{tikzpicture}
\caption{The Dynkin diagrams for the finite Weyl groups. Cominuscule nodes are boxed.}
\label{fig:irreducible-root-systems}
\end{figure}
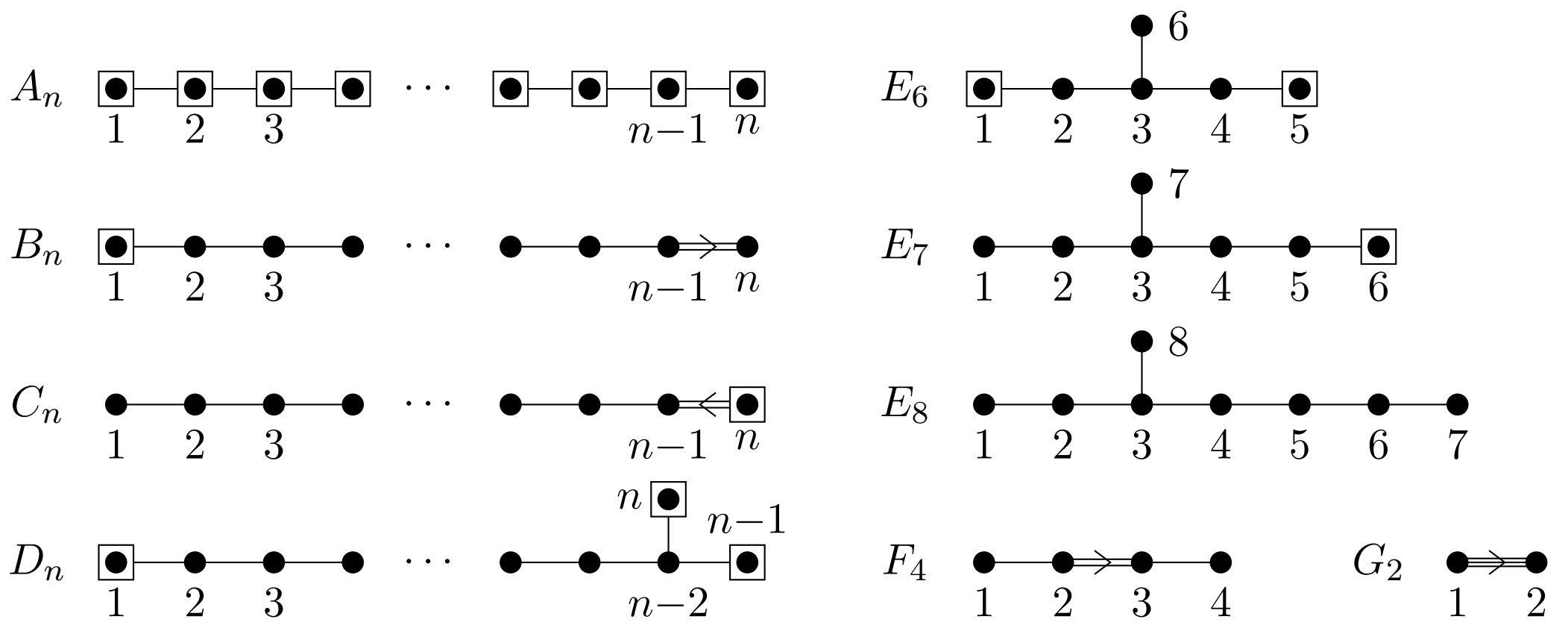

In this section, let $\Phi$ be a finite irreducible crystallographic root system with simple roots $\Delta=\{\alpha_1,\ldots,\alpha_n\}$. For a positive root $\beta=\sum c_i\alpha_i$, write $[\alpha_i]\beta=c_i$ for the coefficient in the basis $\Delta$. Let $\xi$ be the \emph{highest root} of $\Phi$. We say that a simple root $\alpha_i$ is \emph{cominuscule} if $[\alpha_i]\xi=1$. Our conventions for labeling Dynkin diagrams are shown in \Cref{fig:irreducible-root-systems}.

We first prove the forward direction of \Cref{thm:degree-monotone-finite-types}, a strengthening of \cite[Thm.~1.4]{BG}.
\begin{proposition}\label{prop:cominuscule}
Let $J\subseteq S$ be such that $S\setminus J$ contains only cominuscule nodes. Then $W^J$ is degree-monotone.
\end{proposition}
\begin{proof}
Fix $x^J\lessdot y^J\leq w^J$ in $W^J$. Let $x=x^Jw_0(J)$, $y=xs_{\alpha}=y^Jy_J$, and $w=w^Jw_0(J)$. By \Cref{lem:degree-quotient-to-full}, we need to show that $d_w(x)\geq d_w(y)$. By \eqref{eq:plane-sum}, it suffices to show that $d_{M,w}(x)\geq d_{M,w}(y)$, for any $M\in\mathcal{D}_{\alpha}$, a maximal dihedral reflection subgroup containing $s_{\alpha}$. Let $z\in W$ be the minimal element of the coset $xM=yM$, and assume for the sake of contradiction that $d_{M,w}(x)<d_{M,w}(y)$. This cannot happen if $M$ is of type $A_1\times A_1$, so $M$ is either of type $A_2$ or of type $B_2$, as type $G_2$ does not contain any cominuscule nodes. Let the intrinsic generators of $M$ be $s_{\beta}$ and $s_{\gamma}$ with $\beta,\gamma\in\Phi^+$. In both cases, $\beta+\gamma\in\Phi^+$.

We first consider the case where $M$ is of type $A_2$. The only possibility for $d_{M,w}(x)<d_{M,w}(y)$, is that $zM\cap[e,w]=\{z,z_1,z_2,z_3,z_4\}$ where $z< z_i< z_j$ for $i\in\{1,2\}$ and $j\in\{3,4\}$, and that $x=z$ while $y\in\{z_1,z_2\}$. Let $z_5$ be the maximum element of $zM$. For each $\alpha_i\in S\setminus J$, $[\alpha_i](\beta+\gamma)\leq[\alpha_i]\xi=1$ so we can without loss of generality assume that $[\alpha_i]\beta=0$. Either $z_5=z_3s_{\beta}$ or $z_5=z_4s_{\beta}$ so $z_5^{K_i}=\max(z_3^{K_i},z_4^{K_i})$ for $K_i=S\setminus\{\alpha_i\}\subseteq S$. Since $w\geq z_3,z_4$, we have $w^{K_i}\geq z_5^{K_i}$ for each $\alpha_i\in S\setminus J$. Since $\cap _{\alpha_i\in S\setminus J}K_i=J$, we have $w^J\geq z_5^J$ by \cite[Thm.~2.6.1]{BB}. But then $w\geq z_5$, as $w_J=w_0(J)$, a contradiction. 

Lastly assume that $M$ is of type $B_2$, and write $zM=\{z=z_0,z_1,z_2,\ldots,z_7\}$, sorted by length. Here, $W$ is of type $B$ or $C$, both of which contain exactly one cominuscule node $\alpha_i$, with $J=S\setminus\{\alpha_i\}$. Since $[\alpha_i](\beta+\gamma)=1$, we may again assume that $[\alpha_i]\beta=0$ so that $s_{\beta}\in W_J$. Since $x=x^Jw_0(J)$ and $xs_{\beta}<x$, we know that $x\neq z$. The only possibility for $d_{M,w}(x)<d_{M,w}(y)$ is that $zM\cap[e,w]=\{z,z_1,z_2,\ldots,z_6\}$, and that $x\in\{z_1,z_2\}$ and $y\in\{z_3,z_4\}$. Either $z_7=z_5s_{\beta}$ or $z_7=z_6s_{\beta}$, so similar to the above argument, we have $w^J\geq\max(z_5^J,z_6^J)=z_7^J$, and $w\geq z_7$, contradicting $z_7\notin zM\cap[e,w]$. 
\end{proof}

\begin{proof}[Proof of \Cref{thm:degree-monotone-finite-types}]
The cases (1) and (2) are taken care of by \Cref{prop:cominuscule}. Interestingly, type $B_n$ and $C_n$ have isomorphic Coxeter groups, but different cominuscule nodes, allowing us to deduce degree-monotonicity for genuinely different parabolic quotients. For case (3), $(H_3,I_2(5))$ is a finite check. For dihedral Coxeter groups $W$, the parabolic Bruhat graph $\Gamma^J_w$ is always regular. 

Now suppose that $W^J$ is degree-monotone. Write $w_0=w_0^J w_0(J)$ so that $w_0^J$ is the maximal element in $W^J$. In order to narrow down to the cases listed in the theorem statement, we make use of two reductions. First, if $W^J$ is degree-monotone, then $W_K^{J\cap K}$ is degree-monotone for any $K\subseteq S$, since we are just restricting to $w \in W^J$ with $\supp(w) \subseteq K$. Second, if $W^J$ is degree-monotone, then $W^K$ is degree-monotone for $J\subseteq K\subseteq S$ by \Cref{lem:W-K-degree-monotone}. We only need note counterexamples in the remaining cases; these are given in \Cref{tab:monotonicity-counterexamples}. 

\begin{table}[htbp]
\centering
\small
\setlength{\tabcolsep}{3pt}
\renewcommand{\arraystretch}{1.2}
\begin{tabular}{@{}|c|c|c|c|c|c|@{}}
\hline
Type of $W$ & $S\setminus J$ & $w^J$ & $x^J$ & $y^J$
  & $d_{w^J}^J(x^J)<d_{w^J}^J(y^J)$ \\
\hline
$B_n$ & $\{s_1,s_n\}$ & $s_{n-1}s_ns_{n-2}s_{n-1}s_np_{n-2}$
  & $s_np_{n-3}$ & $s_np_{n-2}$ & $n+3<n+4$ \\
$B_n$ & $\{s_2\}$ & $s_2w_0^J$
  & $e$ & $s_2$ & $4n-6<4n-5$ \\
$D_n$ & $\{s_2\}$ & $s_2w_0^J$
  & $e$ & $s_2$ & $4n-8<4n-7$ \\
$E_6$ & $\{s_2\}$ & $s_4s_2s_3s_6s_3s_4s_2s_3s_1s_2$
  & $s_1s_2$ & $s_3s_1s_2$ & $10<11$ \\
$E_6$ & $\{s_6\}$ & $s_6w_0^J$
  & $e$ & $s_6$ & $20<21$ \\
$E_7$ & $\{s_1\}$ & $s_1w_0^J$
  & $e$ & $s_1$ & $32<33$ \\
$E_8$ & $\{s_7\}$ & $s_7w_0^J$
  & $e$ & $s_7$ & $56<57$ \\
$F_4$ & $\{s_1\}$ & $s_1w_0^J$
  & $e$ & $s_1$ & $14<15$ \\
$F_4$ & $\{s_2\}$ & $s_3s_1s_2s_3s_1s_2$
  & $e$ & $s_2$ & $6<7$ \\
$H_3$ & $\{s_1\}$ & $s_1w_0^J$
  & $s_2s_1$ & $s_1s_2s_1$ & $11<12$ \\
$H_3$ & $\{s_2\}$ & $s_2w_0^J$
  & $e$ & $s_2$ & $12<13$ \\
$H_4$ & $\{s_3\}$ & $s_1s_2s_3s_4s_1s_2s_1s_2s_3$
  & $s_4s_3$ & $s_4s_2s_3$ & $10<11$ \\
$H_4$ & $\{s_4\}$ & $s_4w_0^J$
  & $e$ & $s_4$ & $44<45$ \\
\hline
\end{tabular}
\caption{For $B_n$, take $n\geq 3$, and let $p_k=s_ks_{k-1}\cdots s_1$. For $D_n$, take $n\geq 4$.
For $H_3$ and $H_4$ we use the numbering with $(s_1s_2)^5=1$. Computations of the degrees $d^J$ are routine and omitted.}
\label{tab:monotonicity-counterexamples}
\end{table} 
\end{proof}

\section*{Acknowledgments}
C.G. was partially supported by National Science Foundation award no. DMS-2452032 and by a travel grant from the Simons Foundation. Y.G. was partially supported by NSFC Grant no. 12471309. We are also thankful for the excellent collaborative environment provided by the Sydney Mathematical Research Institute, which we visited in May 2025 for the program ``Modern Perspectives in Representation Theory'', and where some of this work was conducted.

\section*{Statement on AI use}
The authors began this project in 2024, so several of the main results were proven without the aid of AI tools and several were AI-assisted. In particular, \Cref{thm:grassmannian}, \Cref{cor:grassmannian-lookup}, and \Cref{thm:degree-monotone-finite-types} were conceived and proven by the authors without the aid of AI.

The authors developed the overall approach to \Cref{thm:main} via maximal dihedral reflection subgroups, Equation (\ref{eq:plane-sum}), and reduction of the problem to studying slices of Bruhat intervals by rank three reflection subgroups. GPT 6 Astra supplied the proofs of several results needed to complete this approach, including the positivity $\delta_M \geq 0$ in \Cref{prop:plane-comparison}, \Cref{prop:separation}, and the main geometric input, \Cref{prop:slice-geometry}.

Claude Opus 5.5 was used for some editing. The authors remain responsible for the content of the manuscript.

\bibliographystyle{amsalpha}
\bibliography{lookup}
\end{document}